\documentclass[a4paper,11pt]{article}
\makeatletter
   \let\accent@spacefactor\relax
 \makeatother       
\usepackage[french]{babel}
\usepackage{amsmath,amscd,amssymb}
\usepackage{fancyheadings}

\def \a {\alpha}
\def \M {M_{{\p}^3}(0,3,0)}
\def \p {{\mathbb P}}
\def \pu {{\p}^1}
\def \noi {\noindent}
\def \vs {\vskip}
\def \bI {{\bf I}}
\def \I {{\cal I}}
\def \di {\partial{\bf I}^1_3}
\def \dj{{\bf \partial I}^2_3}
\def \U {{\bf U}}
\def \Ma {M_{{\p}^3}(0,3,2)}
\def \oo {{\cal O}}
\def \P {{\bf \Psi}}
\def \G {{\mathbb G}}

\begin{document}


\renewcommand{\thesection}{\arabic{section}}

\def \I {{\cal I}}
\def \pd {{\check{\p}^3}}
\def \fl {\longrightarrow}
\def \om {\Omega}
\def \e {\underline{\bf Ext}}
\def \h {\underline{\bf Hom}}
\def \pvv {{\p}({\check V})}
\def \pw {{\p}(W)}
\def \vs {\vskip}
\def \F {{\bf F}}
\def \T {{\bf T}}
\def \dm {{\textit{D\'emonstration}}}
\def \ot {\otimes}
\def \R {{\bf R}}
\def \H {\mathfrak H}
\def \LV {\Lambda^2V}
\def \V {{\cal V}}
\def \g {{\cal G}}
\def \L {{\cal L}}
\def \noi {\noindent}
\def \f {{\bf f}}
\def \cv {{\check V}}
\def \vp {\varphi}
\def \hom {{\rm{Hom}}}
\def \ext {{\rm{Ext}}^1}
\def \oo {{\cal O}}
\def \op {{\cal O}_{\p^3}}
\def \oopu {\oo_{\pu}}
\def \C {\mathbb{C}}

\centerline{\Large{\bf Deux composantes du bord de $\bI_3$}}
\vs 0.1 cm

\centerline{\large{\textsc{Nicolas Perrin}}}
\centerline{{Universit\'e de Versailles}}
\vs -0.1 cm
\centerline{45 avenue des \'Etats-Unis}
\vs -0.1 cm
\centerline{78035 Versailles Cedex}
\vs -0.1 cm
\centerline{email : \texttt{perrin@math.uvsq.fr}}

\vs 0.5 cm

\section*{Introduction}

Un instanton $E$ de degr\'e $n$ est un fibr\'e vectoriel stable de
rang $2$ sur $\p^3$ tel que $c_1(E)=0$, $c_2(E)=n$ et tel que le
groupe $H^1E(-2)$ est nul. La vari\'et\'e $\bI_n$ des instantons de
degr\'e $n$ est donc un ouvert de l'espace des modules ${\bf
M}_{\p^3}(0,n,0)$ des faisceaux sans torsion semi-stables de classes
de Chern $(0,n,0)$. Dans cet article nous nous int\'eressons au bord
de $\bI_n$ dans la vari\'et\'e ${\bf M}_{\p^3}(0,n,0)$.

\vs 0.1 cm

On ne sait d\'ecrire ce bord 
que pour $n=1$ \cite{Ba1} et $n=2$ \cite{NT}. 
%
%
%
%
%
%
La vari\'et\'e $\bI_3$ est de dimension 21. Elle a \'et\'e \'etudi\'ee
par G. Ellingsrud  et S.A. Str\o mme \cite{ES} qui ont prouv\'e
son irr\'eductibilit\'e et par L. Gruson et M. Skiti  \cite{GS} qui
ont montr\'e qu'elle est birationnelle aux r\'eseaux de quadriques
de $\pd$. Au vu des r\'esultats de \cite{NT}, L. Gruson et G. Trautmann ont
conjectur\'e que le bord de $\bI_3$ a huit composantes
irr\'eductibles toutes en codimension 1 (voir remarque 7 pour plus de
d\'etails sur cette conjecture). Dans
\cite{GS}, les auteurs
mettent \'egalement en \'evidence deux composantes irr\'eductibles du
bord de $\bI_3$ correspondant \`a des diviseurs des r\'eseaux de
quadriques. Nous donnons ici une g\'en\'eralisation de ces r\'eseaux
et d\'ecrivons deux nouvelles composantes
irr\'eductibles $\di$ et $\dj$ du bord de la vari\'et\'e des
instantons de degr\'e 3 qui sont ``sym\'etriques'' l'une de l'autre.

\vs 0.1 cm

Consid\'erons la famille $\di$ des faisceaux obtenus comme noyau d'une
fl\`eche surjective de $E''$ vers $\theta(2)$ o\`u $E''$ est un
instanton de degr\'e $1$ et $\theta$ est une th\'eta-caract\'eristique
ayant pour support une conique lisse. Cette famille est contenue dans
$\M$. Consid\'e\-rons par ailleurs un ouvert $\U$ (que l'on d\'efinira
plus loin) de $\Ma$. Un faisceau g\'en\'eral de $\U$ a deux points
singuliers. Soit $\U'$ le ferm\'e  qui correspond aux faisceaux $E''$
de $\U$ dont le lieu singulier est un point double, nous montrons
alors le th\'eor\`eme suivant :

\vs 0.3 cm

\def \th #1#2{\noi\textsc{#1 #2.} ---}

\noi
\th{Th\'eor\`eme}{\!} \textit{(\i) La vari\'et\'e $\di$ est birationnelle \`a
la vari\'et\'e des surfaces quintiques rationelles r\'egl\'ees
autoduales.
C'est une composante irr\'eductible de dimension 20 
du bord de $\bI_3$.}

\vs 0.1 cm

\textit{(\i\i) La famille $\U$ est une fibration en $\p^9$
au dessus du sch\'ema des cubiques gauches irr\'edu\-cti\-bles. Le ferm\'e
$\U'$ est 
d\'ecrit dans chaque fibre par une hypersurface irr\'eductible de
degr\'e $6$.}

\vs 0.1 cm

\textit{(\i\i\i) Il existe une application rationnelle naturelle,
g\'en\'eriquement injective, du ferm\'e $\U'$ dans $\M$ dont l'image
$\dj$ est une composante irr\'eductible du bord de $\bI_3$. L'application
r\'eciproque associe \`a un faisceau $E\in\dj$ son bidual $E''\in\U'$.}

\vs 0.3 cm


Pour prouver que ces familles sont adh\'erentes \`a la vari\'et\'e des
instantons, nous montrons qu'elles peuvent \^etre param\'etr\'ees
par des diviseurs dans une g\'en\'eralisation des r\'eseaux de
quadriques (paragraphe 1). Cette param\'etrisation fait appara\^itre
la sym\'etrie 
entre $\di$ et $\dj$.

\vs 0.1 cm

Dans les deuxi\`eme et troisi\`eme paragraphes, nous \'etudions
respectivement les familles $\di$ et $\dj$. Nous d\'ecrivons en
particulier la g\'eom\'etrie de leurs \'el\'ements de saut.
Rappelons que l'on appelle plan instable (resp. droite bisauteuse) de
$E\in\M$ tout plan $H$ (resp. toute droite $L$) tel que $h^0E_H>0$
(resp. telle que $h^1E_L>0$). 
La courbe des plans instables des faisceaux de ces deux composantes
est ACM de degr\'e 6 et de genre 3 et la surface r\'egl\'ee recouverte
par les droites bisauteuses est soit une quintique rationnelle soit
une sextique elliptique.
%
%
Les \'el\'ements de saut des deux composantes sont reli\'es
ce qui traduit la 
sym\'etrie
entre $\di$ et $\dj$. 
Elle s'exprime assez simplement en termes
de transformations cubo-cubiques (voir paragraphe 1).
En effet, une courbe ACM de degr\'e $6$ et de
genre $3$ permet de d\'efinir 
une application
birationnelle --- appel\'ee transformation cubo-cubique --- de $\p^3$
dans lui m\^eme.
Son application r\'eciproque est encore une
transformation cubo-cubique (i.e. 
associ\'ee \`a une courbe
ACM de degr\'e $6$ et de genre $3$). La 
\textit{dualit\'e}
se traduit par le fait que les transformations cubo-cubiques qui
proviennent de $\dj$ sont les inverses 
de celles 
de $\di$.

\vs 0.1 cm

Dans le troisi\`eme paragraphe, nous \'etudions la famille $\U$
que nous d\'ecrivons explicitement ainsi que le ferm\'e $\U'$ et le
ferm\'e de $\U$ des faisceaux $E''$ qui ne sont plus r\'eflexif. Nous
donnons alors la courbe qui constitue le lieu singulier de $E''$.

\vs 0.1 cm

Enfin dans un dernier paragraphe nous pr\'esentons deux situations
g\'eom\'etriques remarquables reli\'ees \`a notre \'etude. 
Nous donnons une description birationnelle de l'espace des modules des
courbes de degr\'e $7$ et de genre $2$ d'un espace projectif. Cette
description et la formule d'H\" urwitz permettent de retrouver le fait
que $\U'$ est d\'efinie dans les fibres au dessus des cubiques gauches
par une hypersurface de degr\'e $6$.

\vs 0.4 cm

\th{Proposition}{\!} \textit{L'espace des modules des courbes de
degr\'e 7 et de genre 2 de $\p^3$ est birationnellement isomorphe au
quotient par $PGL_2$ de la vari\'et\'e $\G(2,H^0\oo_{S^2\pu}(3))$ des
pinceaux de cubiques du plan $S^2\pu$.}

\vs 0.2 cm


Nous exhibons \'egalement une famille $\mathfrak{I}$ de dimension 36
d'involutions birationnelles de $\p^3$ dans lui-m\^eme. Nous montrons
que $\mathfrak{I}$ est birationnelle au sch\'ema de Hilbert $\H_{9,6}$
des courbes de degr\'e 9 et de genre 6. Nous montrons par ailleurs :

\vs 0.4 cm

\th{Proposition}{\!} \textit{L'espace des modules des courbes de
degr\'e 9 et de genre 6 de $\p^3$ est birationnellement isomorphe au
quotient par $PGL_2$ de la vari\'et\'e $\G(4,H^0\oo_{S^2\pu}(3))$ des
sous-espaces vectoriels de dimension 4 de cubiques du plan $S^2\pu$.}


\vs 0.4 cm

\noi
{\bf Remerciements} : Je tiens \`a remercier ici mon directeur de
th\`ese \textsc{Laurent Gruson} pour toute l'aide qu'il m'a apport\'ee
durant la pr\'eparation de ce travail.

\section{Construction et d\'eformation de faisceaux}

Dans ce paragraphe nous expliquons comment la suite spectrale de
Beilison (voir \cite{OSS}) ramene les probl\`emes sur les
faisceaux \`a des probl\`emes d'alg\`ebre lin\'eaire. Nous \'etudions
alors ces questions d'alg\`ebre lin\'eaires et leurs traductions
g\'eom\'etriques. En particulier nous construisons ainsi des
d\'eformations de faisceaux.

\subsection*{Les transformations cubo-cubiques}
\addcontentsline{toc}{subsection}{Les transformations cubo-cubiques}

\def \X {{\bf X}}

Soit $R$ un espace vectoriel de dimension trois, soient $V$ et $W$ deux
espaces vectoriels de dimension quatre. Une transformation cubo-cubique est
un \'el\'ement $A$ de 
${\bf T}=\p({\rm{Hom}}(R\otimes V, W))$. 
C'est une application lin\'eaire de $R$ dans
$H^0{\cal O}_{{\p}({\check V})\times {\p}(W)}(1)$ dont l'image 
d\'efinit une sous-vari\'et\'e $\Pi$ de ${\p}({\check
V})\times {\p}(W)$. La projection $p$ (resp. $q$) de $\Pi$ sur
${\p}({\check V})$ (resp. ${\p}(W)$) est pour une transformation
cubo-cubique g\'en\'erale
l'\'eclatement de la courbe $Y$ (resp. $Y'$) ACM de degr\'e 6 et de
genre 3 donn\'ee par la r\'esolution :
$$0 \fl R\otimes \oo_{{\p}({\check V})}(-4) \stackrel{A}{\fl}
W\otimes \oo_{{\p}({\check V})}(-3) \fl {\cal I}_Y \fl 0$$
\vs -0.2 cm
respectivement :
\vs -0.4 cm
$$0 \fl R\otimes \oo_{{\p}(W)}(-4) \stackrel{^tA}{\fl} {\check
V}\otimes \oo_{{\p}(W)}(-3) \fl {\cal I}_{Y'} \fl 0$$
Cette construction d\'ecrit une application birationnelle de
${\p}({\check V})$ dans ${\p}(W)$. Si on \'echange les r\^oles de
$\cv$ et de $W$ on a son application r\'eciproque. Nous dirons qu'une transformation cubo-cubique associ\'ee \`a une
courbe $Y$ est involutive si il existe un isomorphisme $\a:W\fl\cv$
tel que $\a(Y')=Y$. La transformation cubo-cubique et son inverse sont
alors d\'efinies par la m\^eme courbe. 
%
%
Pour
plus de d\'etails sur les transformations cubo-cubiques voir aussi
\cite{RS} p. $179$.

\vs 0.4 cm

\th{Remarque}{1} Le groupe $PGL(R)\times PGL(W)$ agit sur $\T$. Il
existe un bon quotient pour cette action qui est fibr\'e principal
homog\`ene sous ce groupe sur un ouvert. Par ailleurs, sur l'ouvert
des \'el\'ements de $\T$ qui d\'efinissent une courbe,  
il existe un morphisme $f$ vers le sch\'ema de Hilbert
$\H_{6,3}$ des courbes ACM de degr\'e 6 et de genre 3 : \`a $A$ on
associe $Y$. Ce morphisme est sur cet ouvert le bon quotient de $T$
par $PGL(R)\times PGL(W)$ (pour plus de d\'etails voir \cite{ESGIT}).


\subsection*{G\'en\'eralisation des r\'eseaux de quadriques}
\addcontentsline{toc}{subsection}{G\'en\'eralisation des r\'eseaux de
quadriques}

Dans la suite on identifie $V$ \`a l'espace vectoriel
$H^0\oo_{\p^3}(1)$. Dans [GS], la famille des instantons sans droite
trisauteuse est identifi\'ee \`a un ouvert des r\'eseaux de quadriques
de $\pd$. Cette identification se fait gr\^ace aux multiplications du
module de Rao d'un faisceau $E\in\bI_3$ : pour un instanton
g\'en\'eral $E$, on a $h^1E(-1)=3$, $h^1E=4$ et $h^1E(1)=1$. La
multiplication $H^1E\ot V\fl H^1E(1)$ est non d\'eg\'en\'er\'ee, elle
permet d'identifier $H^1E$ \`a ${\check V}$. La multiplication 
$$H^1E(-1)\ot V\fl H^1E\backsimeq {\check V}$$
\vs -0.1 cm
\noi
donne alors le r\'eseau de quadriques. 

\vs 0.1 cm

Nous n'allons plus maintenant identifier
directement $H^1E$ \`a ${\check V}$ mais \textit{garder en m\'emoire}
cette identification par la donn\'ee d'un morphisme (qui est un
isomorphisme dans le cas des instantons) de $H^1E$ dans ${\check
V}$ et de sa ``r\'eciproque'' de ${\check V}$ dans $H^1E$. Les
composantes du bord vont appara\^itre lorsque ces morphismes ne
seront plus des isomorphismes. Ainsi pour retrouver $E$ nous aurons besoin
d'identifier $H^1E$ (resp. $H^1E(-1)$) \`a un espace
vectoriel $W$ (resp. $R$) de dimension $4$ (resp.
$3$), d'un morphisme $R\fl W\ot{\check V}$ (une transformation
cubo-cubique) et de deux morphismes $W\fl{\check V}$ et ${\check V}\fl
W$ v\'erifiant les conditions de sym\'etrie suivantes : la
compos\'ee $R\fl {\check V}\ot{\check V}$ (resp. $R\fl W\ot W$) se
factorise par $S^2{\check V}$ (resp. $S^2W$). Dans le cas des
instantons la premi\`ere condition traduit le 
fait que la multiplication dans le module
de Rao est associative. Nous faisons intervenir la seconde
fl\`eche pour conserver la sym\'etrie : sur l'ouvert des instantons,
la fl\`eche de $W$ dans ${\check V}$ est inversible et
l'identification nous permet de faire jouer le m\^eme r\^ole \`a $W$
et ${\check V}$. Nous avons ainsi une fl\`eche de $\cv$ dans $W$
(l'inverse de la pr\'ec\'edente) telle que la compos\'ee $R\fl
W\ot W$ se factorise par $S^2W$. 

Nous nous donnons donc pour g\'en\'eraliser les
r\'eseaux de quadriques trois fl\`eches $R\ot V\stackrel{\vp}{\fl} W$,
$W\stackrel{\psi}{\fl}\cv$ et $\cv\stackrel{\psi'}{\fl} W$ qui
v\'erifient les conditions de sym\'etrie pr\'ec\'edentes.
Consid\'erons ainsi la sous-vari\'et\'e $\F$ de 
${\bf T}\times\p({\rm{Hom}}(W,{\check V}))\times
\p({\rm{Hom}}(\cv,W))$
form\'ee des triplets $(\varphi ,\psi, \psi' )$ tels que les
compos\'ees :
$$R\otimes V\otimes V \stackrel{\varphi\otimes {\bf 1}_V}{\fl}
W\otimes V \stackrel{\psi\otimes {\bf 1}_V}{\fl} {\check V}\otimes V
\fl {\bf C}\ \ \ {\rm{et}} \ \ \ R\otimes {\check W}\otimes {\check W}
\stackrel{\varphi\otimes {\bf 1}_{\check W}}{\fl} \cv\otimes {\check
W} \stackrel{\psi'\otimes {\bf 1}_{\check W}}{\fl} W\ot{\check W} \fl
{\bf C}$$ 
se factorisent par $R\otimes S^2V$ et $R\otimes S^2{\check W}$ et tels
que : 
$$\psi'\circ\psi=\lambda {\bf 1}_{W}\ \ {\rm et }\ \
\psi\circ\psi'=\mu {\bf 1}_{\cv}$$
ces derni\`eres conditions viennent de l'identification au niveau des
instantons : $\psi$ et $\psi'$ sont dans ce cas inverses l'un de
l'autre.
En particulier 
d\`es que l'un de ces deux morphismes n'est plus 
inversible les
compos\'ees doivent \^etre nulles. Nous appellons $\F_{i,j}$ l'image
r\'eciproque dans $\F$ du localement ferm\'e de $\p({\rm{Hom}}(W,{\check
V}))\times\p({\rm{Hom}}(\cv,W))$ form\'e des couples d'applications
lin\'eaires $(\psi,\psi')$ de rangs $i$ et $j$. La derni\`ere
condition impose que si $i\not =4$ alors $j\leq 4-i$. Un point
g\'en\'eral de $\F_{i,j}$ pour $i$ fix\'e diff\'erent de $4$ est 
dans $\F_{i,4-i}$ que nous noterons $\F_i$. Notons $\F_4$
l'ouvert sur lequel $\psi$ est inversible.

\subsection*{Lien avec les faisceaux}
\addcontentsline{toc}{subsection}{Lien avec les faisceaux}

Un \'el\'ement 
de $\F_{i,j}$ nous permet gr\^ace \`a la condition de sym\'etrie sur
$\psi$ de construire un complexe :
$$R\otimes \Omega^2(2) \stackrel{\varphi}{\fl} W\otimes \Omega^1(1)
\stackrel{\psi}{\fl} \oo_{{\p}(V)}$$
Sa cohomologie au centre est est \'el\'ement de $\M$. Cette
construction  
d\'etermine 
une application rationnelle $g$ (d\'efinie l\`a o\`u $\varphi$ est
injective et $\psi$ surjective) de
${\bf F}$ vers $\M$. 

Cette application envoie $\F_4$ dans $\bI_3$ et m\^eme sur l'ouvert des
instantons sans droite trisauteuse (voir \cite{GS}). De plus, deux
diviseurs de $\F_{4}$ font appara\^itre deux composantes
irr\'eductibles du bord de ${\bI}_3$ (cf. \cite{GS}). M. Skiti
\cite{Sk} a \'egalement montr\'e que $\F_2$ s'envoie sur les
instantons \`a droite trisauteuse.
Nous montrons que $\F_3$ et $\F_1$ permettent 
d'identifier deux nouvelles composantes du bord de ${\bI}_3$. 

En
particulier, nous regardons l'image par $g$  de $\F_3$ dans $\M$ qui
sera l'un des deux bords recherch\'es ($\di$). Nous \'etudions les
transformations cubo-cubiques associ\'ees. Notons ${\T_3^1}$ l'image
de $\F_3$ dans $\T$. La vari\'et\'e $\F_4$ correspond au cas o\`u 
la transformation cubo-cubique 
est involutive. Ce n'est plus le cas pour $\F_3$, on a une rupture de
sym\'etrie. Les transformations cubo-cubiques inverses de celles de
$\T_3^1$ sont celles obtenues comme image de $\F_1$ (proposition 12). Elles nous permettrent de d\'ecrire la seconde composante
du bord de $\bI_3$.

Sur $\F_{1}$, la seconde fl\`eche du complexe
n'est plus surjective et nous avons alors une application rationnelle
vers $\Ma$. L'application $g$ n'utilise pas la condition de
sym\'etrie sur $\psi'$. Nous pouvons d\'efinir $g$ sur le localement
ferm\'e $\F^1$ de ${\bf T}\times \p({\rm{Hom}}(W,{\check V}))$ des
paires qui v\'erifient la premi\`ere condition de sym\'etrie et telles
que la seconde fl\`eche est de rang $1$. L'application ainsi d\'efinie 
est \`a valeurs dans un ouvert $\U$ (voir paragraphe 3) de $\Ma$. Nous
verrons qu'il est dominant. L'image de $\F_1$ dans $\F^1$ est un ferm\'e et
son image par $g$ est un ferm\'e $\U'$ (voir paragraphe 3) de
codimension $1$ de $\U$ (l\`a o\`u le lieu singulier du faisceau est
concentr\'e en un point). Nous montrons que la seconde condition de
sym\'etrie (sur $\psi'$) nous permet de 
prolonger ce morphisme vers $\M$. 
Son image d\'ecrit la seconde composante irr\'eductible $\dj$ du bord de
$\bI_3$ recherch\'ee.

\vs 0.4 cm 

\th{Remarques}{2} \textit{(\i)} La fibre g\'en\'erique du morphisme $g$ 
est isomorphe \`a $PGL(R)\times PGL(W)$. En effet, le
faisceau obtenu est \`a cohomologie naturelle. Les multiplications de
$H^1_*E$ nous permettent de retrouver $\vp$ et $\psi$. Sur un ouvert
rencontrant $\F_4$ et $\F_3$, le morphisme $\psi'$ est uniquement
d\'etermin\'e par les deux autres, la fibre est alors isomorphe \`a
$PGL(R)\times PGL(W)$ (cf. \cite{GS} pour $\F_4$ et proposition 2 pour
$\F_3$). 

\textit{(\i\i)} Les vari\'et\'es $\F_1$ et $\F_3$ sont isomorphes (il
suffit d'\'echanger les roles de $V$ et $W$). Par ailleurs, les
composantes $\di$ et $\dj$ du bord de $\bI_3$ sont birationnelles au
quotient de $\F_3$ et $\F_1$ par $PGL(R)\times PGL(W)$. Ainsi les
deux quotients $\di/PGL(V)$ et $\dj/PGL(V)$ sont
birationnels (l'existence de ces deux quotients est assur\'ee au moins
sur un ouvert par un th\'eor\`eme de Rosenlicht \cite{Ro}). Ceci est
une autre forme de la sym\'etrie entre $\di$ et $\dj$.

\vs 0.4 cm

\subsection*{Propri\'et\'es de $\F$}
\addcontentsline{toc}{subsection}{Propri\'et\'es de $\F$}

Nous donnons ici une d\'emonstration de quelques propri\'et\'es de la
vari\'et\'e $\F$, des vari\'et\'es $\F_i$ pour $1\leq i\leq 4$ et de
la vari\'et\'e $\F^1$.

\vs 0.4 cm

\noi
\th{Fait}{1} \textit{Les vari\'et\'es $\F_i$ pour $1\leq i\leq 3$
sont irr\'eductibles de dimension $43$.}

\vs 0.2 cm

\dm : 
Soit $Z_i$ l'image de la projection de $\F_i$ vers le produit
d'espaces projectifs $\p({\rm{Hom}}(W,{\check
V}))\times\p({\rm{Hom}}(\cv,W))$. L'image de la projection de $Z_i$ vers
$\p({\rm{Hom}}(W,{\check V}))$ est la sous-vari\'et\'e
irr\'eductible de dimension $15-(4-i)^2$ des morphismes de rang
$i$. La fibre de cette projection au dessus de $\psi$ est alors donn\'ee
par $\p(({\rm{Ker}\psi{\check )}}\otimes{\rm{Coker}\psi})$. La
vari\'et\'e $Z_i$ est donc irr\'eductible de dimension $14$. Mais
alors la fibre au dessus de $(\psi,\psi')\in Z_i$ est donn\'ee par
$\p(R\otimes(S^2{\rm{Im}}\psi\oplus S^2{\rm{Ker}}\psi\oplus
({\rm{Im}}\psi{\check )}\otimes{\rm{Ker}}\psi))$ ainsi $\F_i$ est
irr\'eductible de dimension $43$.

\vs 0.4 cm

\noi
\th{Fait}{2} \textit{Les vari\'et\'es $\F_i$ pour $1\leq i\leq 3$
sont adh\'erentes \`a $\F_4$. Les vari\'et\'es $\F$ et $\F_4$ sont
irr\'eductibles de dimension $44$.}

\vs 0.2 cm

\dm :
Soit $(\vp,\psi,\psi')$ un \'el\'ement de $\F_i$ assez g\'en\'eral. On
d\'ecompose $\cv$ et $W$ en ${\rm{Im}\psi}\oplus{\rm{Coker}\psi}$ et
${\rm{Im}\psi}\oplus{\rm{Ker}\psi}$. On peut alors supposer que, dans
des bases bien choisies, les applications lin\'eaires $\vp$, $\psi$ et
$\psi'$ s'\'ecrivent sous la forme suivante :
$$\left(\begin{array}{cc}
 A & 0 \\
 C & D \\
\end{array}\right),\ \ 
\left(\begin{array}{cc}
 I & 0 \\
 0 & 0 \\
\end{array}\right) \ \ {\rm{et}}\ \ 
\left(\begin{array}{cc}
 0 & 0 \\
 0 & I \\
\end{array}\right)$$
o\`u $A$ et $D$ sont sym\'etriques. Si on se place sur un anneau de
valuation discr\`ete $A$ d'uniformisante $a$, et que l'on consid\`ere
dans les m\^emes bases, les morphismes $\vp_a$, $\psi_a$ et $\psi'_a$
donn\'es par les matrices :
$$\left(\begin{array}{cc}
 A & a^tC \\
 C & D \\
\end{array}\right),\ \ 
\left(\begin{array}{cc}
 I & 0 \\
 0 & aI \\
\end{array}\right) \ \ {\rm{et}}\ \ 
\left(\begin{array}{cc}
 aI & 0 \\
 0 & I \\
\end{array}\right)$$
alors on voit que pour $a$ inversible le triplet d\'efinit un
\'el\'ement de $\F_4$ alors que sa limite est $(\vp,\psi,\psi')$. La
vari\'et\'e $\F_4$ est donc un ouvert dense de $\F$. Il suffit donc de
montrer  son irr\'eductibilit\'e. Soit $Z_4$ l'image de $\F_4$ dans
$\p({\rm{Hom}}(W,{\check V}))\times\p({\rm{Hom}}(\cv,W))$. Le
morphisme de $Z_4$ dans $\p({\rm{Hom}}(W,{\check V}))$ est un
isomorphisme sur l'ouvert des applications inversibles (sa fibre est
donn\'ee par l'inverse). La vari\'et\'e $Z_4$ est donc irr\'eductible
de dimension $15$. De plus la fibre au dessus d'un point de $Z_4$ est
donn\'ee par $\p(R\otimes S^2V)$ avec l'identification de $W$ et $\cv$
gr\^ace \`a $\psi$. La vari\'et\'e $\F_4$ est donc irr\'eductible de
dimension $44$.

\vs 0.4 cm

\noi
\th{Fait}{3} \textit{Le morphisme de $\F^1$ dans $\T$ est
birationnel sur son image.}

\vs 0.2 cm

\dm :
Soit $\vp\in\T$ un \'el\'ement assez g\'en\'eral dans l'image de
$\F^1$. Il existe donc un morphisme $\psi_0$ v\'erifiant les conditions
de sym\'etrie. Dans des bases de $\cv$ et $W$, on peut \'ecrire $\vp$
et $\psi_0$ sous la forme suivante :
$$\left(\begin{array}{cc}
 A & 0 \\
 C & D \\
\end{array}\right) \ \ {\rm{et}}\ \  
\left(\begin{array}{cc}
 I & 0 \\
 0 & 0 \\
\end{array}\right)$$
o\`u $A$ est de taille $1\times 1$ et $D$ de taille $3\times 3$. On
cherche tous les morphismes $\psi$ de rang $1$ qui v\'erifient les
conditions de sym\'etrie. Un tel morphisme peut s'\'ecrire dans les
m\^emes bases sous la forme :
$$\left(\begin{array}{cc}
 \a & \beta \\
 \gamma & \delta \\
\end{array}\right)$$
on a alors les \'equations $\delta D=^t\!\!\!(\delta D)$ et
$\gamma A+\delta C=^t\!\!\!(\beta D)$. La premi\`ere \'equation est un
syst\`eme lin\'eaire de taille $9\times  9$ et si $\vp$ est assez g\'en\'eral
la seule solution est $\delta=0$. De la m\^eme fa\c con, on voit que l'on a
$\beta=0$ et $\gamma=0$ pour $\vp$ assez g\'en\'eral.

\vs 0.4 cm

\th{Fait}{4} \textit{L'action de $PGL(W)\times PGL(R)$ sur les
vari\'et\'es $\F_i$ pour $1\leq i\leq 4$ est libre sur un
ouvert de chacune de ces vari\'et\'es.}

\vs 0.2 cm

\dm :
Un \'el\'ement $(\vp, \psi,\psi')$ de $\F$ est donn\'e par des matrices :
$$\left(\begin{array}{cc}
 A & 0 \\
 C & D \\
\end{array}\right),\ \ 
\left(\begin{array}{cc}
 I & 0 \\
 0 & 0 \\
\end{array}\right) \ \ {\rm{et}}\ \ 
\left(\begin{array}{cc}
 0 & 0 \\
 0 & I \\
\end{array}\right)$$
o\`u $A$ et $D$ sont sym\'etriques et $A$ est de taille $i\times i$ si
$(\vp,\psi,\psi')\in\F_i$. Le quotient par $PGL(R)$ consiste \`a
prendre un sous espace de dimension 3 des matrices de la premi\`ere
forme. Cette action est g\'en\'eriquement libre. Un \'el\'ement $h$ de
$PGL(W)$ qui est dans le stabilisateur v\'erifie
$h\circ\vp=\lambda\vp$, $\psi\circ h=\mu\psi$ et
$h\circ\psi'=\nu\psi'$. Si $h$ est donn\'e par 
$$\left(\begin{array}{cc}
 \a & \beta \\
 \gamma & \delta \\
\end{array}\right)$$
alors on a $\a=\mu I$, $\beta=0$, $\delta=\nu I$, $\a A=\lambda A$,
$\delta D=\lambda D$ et $\gamma A+\delta C=\lambda C$. Ceci nous donne
: $\lambda=\mu=\nu$, $\a=\lambda I$, $\delta=\lambda I$ et $\gamma
A=0$. Si $\vp$ est assez g\'en\'eral on a donc $\gamma=0$ et $h=\lambda
I$.

\vs 0.4 cm

\noi
{\bf Notations} : (\i) Nous noterons ${\G}$ la grassmannienne des droites
de $\p^3$ et $K$ (resp. $Q$) le sous-fibr\'e (resp. le quotient)
tautologique sur $\G$. Nous noterons $\mathbb{F}(i,j;4)$ les vari\'et\'es
d'incidences des quotients de rang $i$ et de rang $j$ de $V$ et $p$ et
$q$ les projections respectives.

(\i\i) Notons $\pi$ le morphisme de $\F$ vers $\T$.

\pagestyle{headings}
\section{La famille $\di$}

%

Consid\'erons la famille $\di$ contenue dans $\M$ form\'ee
des faisceaux obtenus comme noyaux d'une fl\`eche surjective de $E''$
vers $\theta(2)$ o\`u $E''$ est un instanton de degr\'e 1 et $\theta$
est une th\'eta-caract\'eristique sur une conique lisse $C$. Nous
montrons que cette vari\'et\'e forme une famille irr\'eductible
de dimension 20 qui est adh\'erente \`a ${\bI_3}$. C'est une
composante irr\'eductible du bord de $\bI_3$ dans $\M$.

\subsection*{Une composante du bord}

Notons $b$ le morphisme de $\di$ vers $\bI_1$ qui \`a un faisceau $E$
associe son bidual $E''$. Dans le lemme suivant nous rappelons
les propri\'et\'es de la vari\'et\'e $\bI_1$ des instantons de
degr\'e $1$.

\vs 0.2 cm

\noi
\th{Lemme}{1} \textit{La vari\'et\'e $\bI_1$ est isomorphe
\`a ${\p}(\Lambda^2{\check V})$ priv\'e de la grassmannienne ${\G}$ :
la donn\'ee d'un instanton $E''$ de degr\'e 1 est \'equivalente
\`a celle du complexe lin\'eaire non singulier ${\cal A}$ des droites
de saut de $E''$.}

\vs 0.2 cm

\dm :
Nous nous contenterons de d\'ecrire la situation g\'eom\'etrique, pour
plus de d\'etails et une d\'emonstration voir \cite{OSS} page 364. La
vari\'et\'e ${\cal A}$ des droites de saut de $E''$ est un complexe
lin\'eaire 
non singulier de droites (i.e. la trace dans $\G$ d'un hyperplan de
$\p(\LV)$). Il d\'etermine $E''$. Notons $\X$ l'image r\'eciproque par
$q$ de ${\cal A}$ dans $\mathbb{F}(1,2 ;4)$, la vari\'et\'e $\X$
s'identifie \`a $\p_{\p^3}(E''(1))$ et \`a $\p_{\cal A}(Q(1))$.

\vskip 0.4 cm

\noi
\th{Proposition}{1} \textit{La fibre de $b$ au dessus d'un
complexe de droites ${\cal A}$ est birationnellemement isomorphe au
sch\'ema de Hilbert des courbes rationnelles lisses de degr\'e $5$ de
${\cal A}$.}

\textit{La vari\'et\'e $\di$ est donc birationnellemement
isomorphe \`a la vari\'et\'e ${\bf H}^5_5$ des courbes rationnelles lisses
de degr\'e 5 de $\G$ qui sont trac\'ees sur un complexe non singulier
de droites.}

\vskip 0.2 cm

\dm :
Soit $E\in\di$ noyau d'une surjection $E'' \fl \theta(2)$. Notons $C$
le support de $\theta$. La surface r\'egl\'ee $S={\p}_C(E''(1)\vert_C)$
est incluse dans $\mathbb{F}(1,2 ;4)$ et m\^eme dans ${\bf
X}=\p_{\p^3}(E''(1))$. 
La surjection
$E''(1)\vert_C\fl \theta(3)$ d\'efinit une section $\sigma :C\fl S$ de cette
surface. Notons $Z=\sigma(C)$. Pour un faisceau g\'en\'eral, la
projection de $Z$ dans ${\cal A}$ est une immersion ferm\'ee d'image
une courbe $\Gamma$ rationnelle lisse de degr\'e
5. De plus, le $\oopu$-module 
$Q(1)\vert_{\Gamma}$ est isomorphe \`a $\oo_{\pu}(2)\oplus
\oo_{\pu}(3)$. En effet, il est de degr\'e 5 et la courbe $Z$ nous
donne une surjection de $Q(1)\vert_{\Gamma}$ vers $\oo_{\pu}(2)$ (car la
projection par $p$ de $Z$ redonne la conique $C$). La courbe $Z$ est
la section de la surface r\'egl\'ee
${\p}_{\Gamma}(Q(1)\vert_{\Gamma})$ d\'efinie par la fl\`eche
surjective $Q(1)\vert_{\Gamma}\fl\oopu(2)$.

R\'eciproquement, soit $\Gamma$ une courbe rationnelle lisse de
degr\'e $5$ dans ${\cal A}$. Supposons que $\Gamma$ est sur un unique
complexe de droites et que $Q(1)\vert_{\Gamma}$ est isomorphe \`a
$\oo_{\pu}(2)\oplus \oo_{\pu}(3)$ en tant que $\oopu$-module (c'est le
cas g\'en\'erique pour $\Gamma$). Consid\'erons alors la surface
r\'egl\'ee (de morphisme la restriction de $q$) suivante :
${\p}_{\Gamma}(Q(1)\vert_{\Gamma})$, elle est incluse dans ${\bf
X}$. La section donn\'ee par la surjection de $Q\vert_{\Gamma}(1)$
vers le facteur $\oo_{\pu}(2)$ a pour image une courbe 
$Z$ contenue dans ${\bf X}$. Sa projection par $p$ est une conique
$C$ (lisse en g\'en\'eral). Mais alors la courbe $Z$ est une section
de la surface r\'egl\'ee ${\p}E''(1)\vert_C$ (car la projection $p$
est un isomorphisme de $Z$ sur $C$). Cette section correspond \`a une
surjection de $E''(1)\vert_C$ vers un $\oo_{\pu}$-module inversible de
degr\'e celui de $\Gamma$. On a donc une surjection de $E''\vert_C$
vers $\theta(2)$ ce qui nous d\'efinit l'application r\'eciproque.

En prenant le noyau de la
compos\'ee $E'' \fl E''\vert_C \fl \theta(2)$, on retrouve le faisceau $E$.

\vskip 0.4 cm

\noi
\th{Lemme}{2} \textit{La vari\'et\'e ${\bf H}^5_5$ est
irr\'eductible de dimension 20.}

\vskip 0.2 cm

\dm : Cette vari\'et\'e correspond aux surfaces quintiques
rationnelles r\'egl\'ees autoduales.
Elle est irr\'eductibilit\'e de dimension 20 (voir [P2] ou [P3]).

\vskip 0.4 cm 

\noi
\th{Corollaire}{1} \textit{La vari\'et\'e $\di$ est irr\'eductible
de dimension 20.}

\vskip 0.4 cm

Nous montrons gr\^ace aux r\'eseaux de quadriques que la
famille $\di$ est adh\'erente \`a $\bI_3$.

\vs 0.2 cm

\noi
\th{Proposition}{2} \textit{La restriction de $g$ \`a $\F_3$ est
dominante sur $\di$.}

\vskip 0.2 cm

\dm :
Soit $E\in\di$ un faisceau \`a cohomologie naturelle. Nous identifions
$H^1E$ \`a $W$ et $H^1E(-1)$ \`a $R$. La premi\`ere multiplication du
module de Rao nous donne une transformation cubo-cubique
$\varphi\in\T$. Par ailleurs, la seconde multiplication du module de
Rao nous donne un \'el\'ement $\psi\in\p({\rm{Hom}}(W,{\check V}))$
qui v\'erifie la condition de sym\'etrie. De plus nous savons que
$H^1E=H^0\theta(2)$ et $H^1E(1)={\rm{Coker}}(H^0E''(1)\fl
H^0\theta(3))$ donc le morphisme $\psi$ est de rang $3$ (la
multiplication par l'\'equation du plan de $C$ est nulle). Prenons
alors pour $\psi'$ un morphisme de rang $1$ dont le noyau est
${\rm{Im}}\psi$ et l'image est ${\rm{Ker}}\psi$ (ceci revient \`a
choisir un isomorphisme entre ${\rm{Coker}}\psi$ et ${\rm{Ker}}\psi$
qui sont de dimension $1$, il n'y a qu'une solution \`a scalaire
pr\`es). Ce morphisme v\'erifie les conditions de sym\'etrie et
d'annulation. L'\'el\'ement $(\varphi,\psi,\psi')$ est donc dans
$\F_3$. Ainsi $g^{-1}(\di)$ est contenu dans $\F_3$ qui est
irr\'eductible de dimension \'egale \`a ${\rm{dim}}(g^{-1}(\di))$ (la
fibre g\'en\'erale est $PGL(R)\times PGL(W)$ l'action sur $\F_3$
\'etant g\'en\'eriquement libre).
%

\vs 0.4 cm

\noi
\th{Corollaire}{2} \textit{La vari\'et\'e $\di$ est une composante
irr\'eductible du bord de ${\bI}_3$.}

\vskip 0.2 cm

\dm :
Nous savons [GS] que la restriction de $g$ \`a $\F_4$ est une application
rationnelle dominante sur ${\bI}_3$. Or nous avons vu au fait 2 que
$\F_3$ est adh\'erente \`a $\F_4$. Son image par $g$, qui est $\di$,
est donc adh\'erente \`a ${\bI_3}$. Dans \cite{P1} nous d\'ecrivons
$\di$ comme le diviseur exceptionnel d'un \'eclatement.

\subsection*{\'El\'ements de saut}

Soit $E\in\di$, nous \'etudions maintenant les \'el\'ements de saut de
$E$. 

\vs 0.4 cm

\noi
\th{Proposition}{3} \textit{La vari\'et\'e $Y$ des plans instables
de $E$ est une courbe ACM de degr\'e 6 et de genre 3 de $\pd$. Elle a
un point triple au point correspondant au plan de la conique $C$.}

\vskip 0.2 cm

\dm :
La vari\'et\'e d'incidence $\mathbb{F}(1,3;4)$ nous permet de calculer la
r\'esolution de l'id\'eal de la courbe $Y$ des plans instables qui est
donn\'e par $R^1q_*p^*E$ :
$$0 \fl R\otimes{\cal O}_{\pd}(-4) \fl W\otimes {\cal O}_{\pd}(-3) \fl
{\cal I}_Y \fl 0$$
o\`u $R=H^1E(-1)$, $W=H^1E$ et la multiplication est donn\'ee par la
multiplication du module de Rao de $E$. La courbe $Y$ est ACM
de degr\'e $6$ et de genre $3$.

De plus, les groupes de cohomologie de $E$ sont d\'efinis de la fa\c
con suivante :
$$0\fl H^0\theta(1) \fl H^1E(-1) \fl H^1E''(-1) \fl 0\ \ \ {\rm{et}} \
\ \ H^1E=H^0\theta(2)$$ 
Pour voir que $Y$ a un point triple donn\'e par le plan contenant la
conique, nous regardons la matrice $A$ de la multiplication du module de
Rao qui d\'efinit $Y$. La multiplication par l'\'equation $H$ du plan
de la conique a un noyau de dimension 2 dans $H^1E(-1)$. Les mineurs
$2\times 2$ de la matrice $A$ sont donc contenus dans le noyau de $H$
(vu comme forme lin\'eaire sur ${\check V}$). Ainsi, la courbe
$Y$ a un point triple au point de $\pd$ dont l'id\'eal est engendr\'e
par le noyau de $H$.

\vs 0.4 cm

Si $E\in\di$, notons $\Gamma$ la courbe rationnelle quintique
trac\'ee sur $\G$ d\'efinie \`a la proposition 1.

\vs 0.2 cm

\noi
\th{Proposition}{4} \textit{La vari\'et\'e des droites bisauteuses
de $E$ s'identifie \`a celle des tris\'ecantes \`a $Y$ la courbe des
plans instables. C'est une courbe de degr\'e $8$ r\'eunion de $\Gamma$
et d'une cubique du plan dual de celui de la conique $C$.}

\textit{La courbe
$Y$ est le lieu double de la surface r\'egl\'ee d\'efinie par $\Gamma$
et le lieu triple de celle d\'efinie par la courbe des droites
bisauteuses.}

\vskip 0.2 cm

\dm :
Une droite $L$ qui passe par le point triple de $Y$ est tris\'ecante
\`a $Y$ (c'est \`a dire v\'erifie la condition $h^1\I_Y(1)\vert_L>0$) 
si elle recoupe $Y$ ou est tangente
\`a $Y$ au point triple. Les tris\'ecantes \`a $Y$ passant par le
point triple forment donc une cubique du plan dual de celui de la conique qui
correspond \`a la projection de $Y$ par rapport \`a son point
triple. De fa\c con g\'en\'erale, les tris\'ecantes \`a $Y$ sont
donn\'ees par le $0^{\rm i\grave{e}me}$ id\'eal de Fitting $Q$ de
$R^1q_*p^*({\cal
I}_Y(1))(1)$ (avec la vari\'et\'e d'incidence $\mathbb{F}(2,3;4)$) qui a la
r\'esolution suivante :
$$R\otimes K \fl W\otimes {\cal O}_{\bf G} \fl Q \fl 0$$

Par ailleurs, les bisauteuses de $E$ sont donn\'ees (gr\^ace \`a la
vari\'et\'e d'incidence $\mathbb{F}(1,2;4)$) par le
$0^{\rm i\grave{e}me}$ id\'eal de Fitting de $R^1q_*p^*E$ qui a la
m\^eme r\'esolution. Les droites bisauteuses sont donc les
tris\'ecantes \`a la courbe $Y$ des plans instables. 

Les droites bisauteuses sont donn\'ees par le support du conoyau de la
fl\`eche suivante : $q_*p^*E''\fl q_*p^*\theta(2)$. Par d\'efinition
de $E''$, en dehors du plan des droites coupant la conique en deux
points, ce support est exactement l'image de la section de $\p_C(E'')$
d\'efinie par la surjection $E''\fl\theta(2)$.
La courbe $\Gamma$ est donc le lieu des tris\'ecantes \`a $Y$ qui ne
passent pas par le point triple. La seconde composante du lieu des
droites bisauteuses est alors donn\'ee par les tris\'ecantes qui
passent par le point triple, c'est la projection de $Y$ \`a partir de
son point triple. C'est une cubique du plan des droites du plan de
$C$.

Il reste \`a d\'eterminer le lieu singulier des surfaces
r\'egl\'ees. Or si $H$ est un plan stable, nous avons la suite exacte
$$0\fl\oo_H\fl E_H(1)\fl\I_Z(2)\fl 0$$ 
o\`u $Z$ est de longueur $4$. Une
droite de $H$ est bisauteuse si elle passe par trois points de $Z$. Il
en existe donc au plus une et les plans stables ne sont pas dans le
lieu singulier. Si $H$ est instable nous avons la suite exacte 
$$0\fl\oo_H\fl E_H\fl\I_Z\fl 0$$
avec $Z$ de longueur $3$. Les trois droites passant par deux des
points de $Z$ sont alors bisauteuses. Les plans instables forment donc
le lieu triple de la surface r\'egl\'ee des bisauteuses. De plus, une
de ces droites est contenue dans le plan de $C$ (car deux des trois
points de $Z$ sont sur la conique) et les deux autres sont sur la
surface r\'egl\'ee qui correspond \`a $\Gamma$.
 
\vs 0.4 cm

\th{Remarque}{3} Notons $\H^t_{6,3}$ (resp. $\H^1_{6,3}$) le ferm\'e
de $\H_{6.3}$ (resp. de $\H^t_{6,3}$) des courbes ayant un point
triple (resp. dont la courbe des tris\'ecantes est trac\'ee sur un
complexe lin\'eaire de droites). 
Les 
vari\'et\'es $\di$, ${\bf H}^5_5$ et $\H^1_{6,3}$ sont birationnelles
et irr\'eductibles. Les morphismes de $\di$ vers ${\bf H}^5_5$ et
$\H^1_{6,3}$ sont donn\'es respectivement par les droites bisauteuses
et les plans instables d'un faisceau. 
Les morphismes entre ${\bf H}^5_5$ et $\H^1_{6,3}$ sont
donn\'es dans un sens par le lieu double de la surface et dans l'autre
par le lieu des tris\'ecantes ne passant pas par le point triple. Cette
derni\`ere correspondance birationnelle est d\'ecrite plus en
d\'etails dans \cite{P2}.




\section{La famille $\dj$}

Nous \'etudions un ouvert ${\bf U}$
de l'espace $M_{{\p^3}}(0,3,2)$ qui param\'etrise les faisceaux de rang 2
sans torsion semi-stables et de classes de chern $c_1=0$, $c_2=3$ et
$c_3=2$ de $\p^3$. Nous montrons que la donn\'ee d'un faisceau $E''$ dans
${\bf U}$ est \'equivalente \`a celle de la surface form\'ee par la
r\'eunion de ses droites bisauteuses et nous d\'ecrivons
g\'eom\'etriquement les surfaces ainsi obtenues.

\vskip 0.4 cm

\noi
\th{Exemple}{1} Soit $Y$ une courbe de degr\'e 7 et de genre 2 assez
g\'en\'erale (notamment lisse et irr\'eductible) et $\xi$ un
\'el\'ement non nul de ${\rm{Ext}}^1_{{\cal O}_{\p^3}}({\cal I}_Y(4),
{\cal O}_{\p^3})$. Il d\'efinit une extension dont le terme central
est un faisceau $E''(2)$. Le faisceau $E''$ est r\'eflexif (voir
\cite{Ha2}) et $E''\in M_{{\p^3}}(0,3,2)$. De plus pour $Y$
assez g\'en\'erale, ces faisceaux sont \`a cohomologie minimale en
degr\'es variant de -2 \`a 1. Ils forment donc un ouvert dans
${\Ma}$.
\vskip 0.4 cm

Soit $\U_0$ l'ouvert de $\Ma$ (qui contient les faisceaux de
l'exemple) des faisceaux $E''$ tels que la cohomologie de $E''(-2)$,
$E''(-1)$, $E''$ et $E''(1)$ est naturelle, c'est \`a dire qu'au plus
un groupe de cohomologie est non nul. On a alors $h^2E''(-2)=1$,
$h^1E(-1)=2$, $h^1E=3$ et tous les groupes de cohomologie de $E(1)$
sont nuls. Cet ouvert ${\bf U}_0$ contient des faisceaux r\'eflexifs
et des faisceaux sans torsion. Nous d\'ecrirons les
sous-vari\'et\'es correspondantes et le lieu singulier des faisceaux
sans torsion. Soit $E''$ dans ${\bf U}_0$ et consid\'erons la suite
spectrale de Beilinson associ\'ee \`a $E''(1)$ (voir par exemple
\cite{OSS}) : $E_1^{p,q}=H^qE''(p+1)\otimes\Omega^{-p}(p) \Rightarrow
E''(1)$. Nous en d\'eduisons les deux suites exactes  :
$$0\fl H^1E''(-1)\otimes\Omega^2(2) \fl H^1E''\otimes \Omega^1(1) \fl F
\fl 0 \ \ \ (1)$$
$$0\fl H^2E''(-2)\otimes \oo_{\p^3}(-1) \fl F\fl E''(1) \fl 0$$

\vskip 0.4 cm

\noi
\th{Remarques}{4} \textit{(\i)} Soit ${\bf W}$ la vari\'et\'e des faisceaux
$F$ obtenus gr\^ace \`a une suite exacte du type de $(1)$. Nous
montrons (proposition 5) qu'elle est birationnelle au sch\'ema de
Hilbert $\H_{3,0}$ des cubiques gauches
irr\'eductibles. La fibre du morphisme de $\U_0$ vers ${\bf W}$ au
dessus de $F$ est un ouvert de $\p(H^0F(1{\check )})$ (l\`a o\`u la
section est injective).
En effet, si $E''(2)$ est le conoyau d'une section de $F(1)$,
alors $E''$ a la cohomologie souhait\'ee.

\textit{(\i\i)} Si $E''$ est r\'eflexif, alors le faisceau $F$ est localement
libre. En effet, soit $Z$, le
lieu d'annulation de la section de $F(1)$, nous avons la suite exacte
$$0\fl \oo_Z(1) \fl \e^1(E''(1),\oo_{\p^3}) \fl \e^1(F,\oo_{\p^3}) \fl 0$$
le sch\'ema $Z$ est n\'ecessairement de codimension 3 et dans ce cas
sa longueur est donn\'ee par $c_3(F)$ qui vaut ici 2. Or
$\e^1(E''(1),\oo_{\p^3})$ est aussi de longueur 2 (sa longueur est
$c_3(E'')$ voir \cite{Ha2}) et nous voyons que la
premi\`ere fl\`eche est un isomorphisme donc le dernier terme de la
suite exacte est nul et $F$ est localement libre.

\subsection*{Etude des \'el\'ements de saut}

Dans la proposition suivante, le faisceau $F$ est donn\'e par un
faisceau $E''\in {\bf U}_0$. 

\vskip 0.4 cm

\noi
\th{Proposition}{5} \textit{La vari\'et\'e des plans instables
de $E''$ est une cubique gauche $C$ de $\pd$. Si $C$ est
irr\'eductible, le dual de $F$ est un fibr\'e de Schwarzenberger
associ\'e \`a $C$.}

\vskip 0.2 cm

\textit{D\'emonstration} :
Les plans instables de $E''$ sont \'evidement ceux de $F(-1)$. Les
formules de Bott et la suite exacte $(1)$ montrent que pour tout plan
$H$ le groupe $H^2F_H(-1)$ est nul. Les
plans instables sont donc les plans $H$ qui v\'erifient la condition
$h^1F_H(-1)>1$. La vari\'et\'e des plans instables est donn\'ee par le
premier id\'eal de Fitting du faisceau
$R^1q_*p^*F_{\mathbb{F}(1,3;4)}(0,-1)$ qui admet la r\'esolution
suivante :
$$0\fl H^1F(-2)\otimes\oo_{\pd}(-1)\stackrel{M}{\fl}
H^1F(-1)\otimes\oo_{\pd}\fl R^1q_*p^*F_{\mathbb{F}(1,3;4)}(0,-1)\fl 0$$
Mais on a $H^1F(-2)=H^1E(-1)$ et $H^1F(-1)=H^1E$ qui sont de dimensions
respectives $2$ et $3$. Ce premier id\'eal de Fitting est une cubique
gauche sauf si les mineurs
de $M$ ont un facteur commun. C'est alors un plan ou une quadrique. Dans
ce cas la fl\`eche de ${\Omega^2(2)}^2$ dans
${\Omega^1(1)}^3$ qui d\'efinit $F$ n'est plus injective ce qui ne
peut se produire. Le faisceau $F$ est compl\`etement d\'etermin\'e par
cette courbe (elle d\'etermine la suite exacte $(1)$).


Six cas sont possibles si la cubique n'est pas irr\'eductible (en
excluant les matrices qui donnent des surfaces) : la courbe est
r\'eunion d'une
conique et d'une droite se coupant en un point ou une chaine de trois
droites ou trois droites concourantes non coplanaires ou une droite
double et une droite non contenue dans son plan qui la recoupe ou une
droite triple dans un c\^one quadratique ou enfin une droite triple
donn\'ee par le carr\'e de l'id\'eal d'une droite (pour la description
de ces courbes, voir \cite{R}). Le faisceau $F$
g\'en\'eral obtenu dans le cas o\`u $C$ d\'eg\'en\`ere est alors sans
torsion de lieu singulier contenant une droite. L'espace vectoriel
$H^0F(1)$ est toujours de dimension $10$ et les faisceaux $E''$
obtenus sont sans torsion et leur lieu singulier contient une
droite. Ces faisceaux forment un ferm\'e de codimension $1$ de
$\U_0$.

Nous nous pla\c cons dans toute la suite sur l'ouvert $\U$ de $\U_0$
form\'e des faisceaux dont la courbe des plans instables est une cubique
gauche irr\'eductible. Cette cubique est isomorphe \`a $\pu$. Notons
alors $S_n$ la repr\'esentation irr\'eductible de dimension $n+1$ de
$SL_2$. 
Dans ce cas, le faisceau ${\check F}(1)$ est un fibr\'e de Steiner :
on a une r\'esolution 
$$0 \fl {\C}^2\ot\oo_{\p^3}(-1) \stackrel{u}{\fl} {\C}^5\ot\oo_{\p^3} \fl
{\check F}(1) \fl 0$$
et les plans instables de $F(-1)$ sont donn\'es par $C$. Les
r\'esultats de \cite{V} nous permettent de dire que ${\check F}(1)$ est un
fibr\'e de Schwarzenberger associ\'e \`a $C$. Les espaces vectoriels
${\C}^2$ et ${\C}^5$ s'identifient alors \`a $S_1$ et $S_4$.

\vskip 0.4 cm

Rappelons que l'on se place maintenant sur l'ouvert $\U$ de $\U_0$
form\'e des faisceaux dont la courbe des plans instables est une
cubique gauche irr\'eductible. Nous \'etudions maintenant la
fibre du morphisme de $\U$ dans ${\bf W}$ au dessus de $F$, c'est \`a
dire aux sections de $F(1)$. 


\vskip 0.4 cm 

\noi
\th{Remarques}{5} \textit{(\i)} L'espace vectoriel $H^0F(1)$ est
isomorphe \`a $S^2S_3$ en tant que $SL_2$-module. L'espace vectoriel
${\check V}$ est isomorphe \`a $S_3$. Ainsi, nous avons
l'identification $\Lambda^2 V=\Lambda^2 S_3=S^2S_2$ (voir
\cite{FH}). La vari\'et\'e des bis\'ecantes 
\`a $C$ d\'ecrit le plongement de Veronese $v$ de $\p^2=\p(S_2)$
dans $\G\subset\p(\Lambda^2V)$ d'image $\V$. 

Dans $\p(S_2)$ nous avons une
conique canonique $C_0$ dont l'image $v(C_0)$ dans $\V$ est la courbe
des tangentes \`a $C$. La donn\'ee de $s\in S^2S_3=S^3S_2$ (loi de
r\'eciprocit\'e de Hermite, voir [FH]) correspond \`a la donn\'ee
d'une cubique $X$ de $\p(S_2)$. La courbe $v(X)$ est alors
elliptique de degr\'e $6$ sur $\V$. Cette courbe d\'efinit une surface
r\'egl\'ee sextique elliplique $S$
dont le mod\`ele non singulier est donn\'e par la restriction du
quotient tautologique $Q$ de la grassmanienne \`a $v(X)$. Le
mod\`ele non singulier de la surface duale ${\check S}$ est donn\'e
par la restriction de ${\check K}$ \`a $v(X)$. 


\textit{(\i\i)} Nous d\'eterminons (proposition 10)
les conditions
n\'ecessaires et suffisantes
sur l'\'el\'ement de $S^2S_3$ (la section de $F(-1)$) pour que le
faisceau $E''$ soit r\'eflexif. 

\vskip 0.4 cm

\noi
\th{Proposition}{6} \textit{La vari\'et\'e des droites de saut de
$F(-1)$ s'identifie \`a la V\'eronese $\V$ des bis\'ecantes \`a $C$.}

\vskip 0.2 cm 

\dm :
La vari\'et\'e d'incidence permet de montrer que les droites sauteuses
sont donn\'ees par le support du faisceau
$R^1q_*p^*F_{\mathbb{F}(1,2;4)}(0,-1)$ dont on a une pr\'esentation :
$$ H^1F(-2)\otimes K \stackrel{\alpha}{\fl} H^1F(-1)\otimes\oo_{\bf G}
\fl R^1q_*p^*F_{\mathbb{F}(1,2;4)}(0,-1)\fl 0$$
Par ailleurs, les bis\'ecantes \`a $C$ sont donn\'ees (cf. \cite{GP}) par
le support du faisceau $R^1q_*J$ o\`u $J$ est l'id\'eal de $p^{-1}C$
dans $\mathbb{F}(2,3;4)$. Or ce faisceau admet la r\'esolution suivante : $S_1\otimes\oo_{\mathbb{F}(2,3;4)}(-3) \fl S_2\otimes
\oo_{\mathbb{F}(2,3;4)}(-2) \fl J \fl 0$ qui nous donne la 
pr\'esentation 
$$S_1\otimes K(-1) \stackrel{\beta}{\fl} S_2\otimes \oo_{\bf
G}(-1) \fl R^1q_*J \fl 0$$
Les fl\`eches $\beta(1)$ et $\alpha$ sont \'egales et les deux faisceaux ont donc le m\^eme id\'eal de Fitting.

Cette \'etude nous permet d'identifier le faisceau $R^1q_*p^*F_{{\bf
F}(1,2;4)}(0,-1)\ot\oo_{\G}(1)$ au faisceau
$R^1q_*J\ot\oo_{\G}(2)$. Ce dernier faisceau est localement libre de
rang $1$ sur $\V$. En effet, si la droite $L$ coupe
$C$ en $a$ points on a $J_L=\oo_L(-a)$. Or une cubique gauche a au plus
des bis\'ecantes 
ainsi $H^1J_L$ est non nul si et seulement si $L$ est bis\'ecante \`a
$C$ et alors $h^1J_L=1$. Ainsi le faisceau $R^1q_*J\ot\oo_{\G}(2)$ est
inversible sur $\V$. La pr\'esentation de $R^1q_*J$ nous permet de
dire que $R^1q_*J\ot\oo_{\G}(2)$ s'identifie \`a
$\oo_{\p(S_2)}(3)$.

\vskip 0.4 cm

\noi
\th{Proposition}{7} \textit{La courbe $v(X)$ est la courbe des
droites bisauteuses de $E''$.}

\vskip 0.2 cm 

\dm :
Notons $\g$ le faisceau
$R^1q_*p^*F_{\mathbb{F}(1,2;4)}(0,-1)\ot\oo_{\G}(1)$.  Nous avons vu
\`a la proposition
pr\'ec\'edente que ce faisceau $\g$ est isomorphe \`a
$\oo_{\p(S_2)}(3)$ d'o\`u l'identification de $H^0\g$ \`a $S^3S_2$. Les
sections de $\g$ et
de $F(1)$ sont les m\^emes et la section $s$ de $F(1)$ nous donne une
section $s_0$ de $\g$. La cubique plane $X$ d\'efinie par $s$ est
exactement le lieu des z\'eros de $s_0$. Or nous avons la suite exacte :
$$\oo_{\bf G} \stackrel{s_0}{\fl} \g \fl
R^1q_*p^*E''_{\mathbb{F}(1,2;4)}\otimes\oo_{\bf G}(1) \fl 0$$
ce qui montre que les bisauteuses de $E''$ sont exactement donn\'ees
par $v(X)$.

\subsection*{Etude de la surface r\'egl\'ee}

Pour aborder la suite de l'\'etude de $\U$, trois points de vue sont
possibles : \'etudier les diff\'erents faisceaux $E''$ obtenus \`a
partir d'une cubique gauche $C$ et d'une cubique plane $X$ de
$\p(S_2)$ ou d\'ecrire la position d'une cubique plane $X$ par rapport
\`a une conique fix\'ee (la conique canonique $C_0$ de $\p(S_2)$) ou
encore d\'ecrire la surface $S$ (ou ${\check S}$) et notamment
pr\'eciser son mod\`ele non singulier.

Les questions qui se posent sur $E''$ sont alors de savoir si $E''$
est r\'eflexif ou non et de d\'ecrire son lieu singulier. Pour la
position de $X$ par rapport \`a $C_0$, dans \cite{P2} nous avons
montr\'e que l'on a
soit deux, soit un unique soit une infinit\'e de triplets sur $X$
associ\'es \`a la conique (sommets d'un triangle de Poncelet tangent
\`a $C_0$). Enfin la surface ${\check S}$ \'etant r\'egl\'ee sur une
courbe elliptique $X$, les r\'esultats de 
[Ha1] (th\'eor\`eme V.2.15) vont nous permettre de dire que le
faisceau ${\check K}_X$, 
localement libre de rang $2$ sur $X$, est soit somme directe de deux
faisceaux inversibles de degr\'es $3$ non isomorphes, soit extension
non triviale de deux faisceaux inversibles de degr\'es $3$ isomorphes,
soit somme directe de deux faisceaux inversibles de degr\'es $3$
isomorphes.

Nous montrons que ces trois probl\`emes sont \'equivalents et que les
diff\'erents cas se correspondent.

\vs 0.4 cm

\noi
\th{Proposition}{8} \textit{Le faisceau ${\check K}_X$ est de l'une des deux formes suivantes :}

\textit{(\i) somme directe de deux faisceaux inversibles de degr\'e $3$}

\textit{(\i\i) extension non triviale d'un faisceau inversible de degr\'e $3$ par lui m\^eme}

\vs 0.2 cm

\dm :
La restriction de la suite exacte tautologique de $\G$ \`a $\V$ nous donne :
$$0\fl S_1\ot\oo_{\p(S_2)}(-1)\fl S_3\ot\oo_{\p(S_2)}\fl {\check
K}_{\V}\fl 0$$
Le fibr\'e ${\check K}_{\V}$ est donc un fibr\'e de Schwarzenberger
\cite{S}.

Le groupe $H^0{\check K}_{X}(-1)$ est nul : il s'identifie \`a
$H^1{\check K}_{\V}(-4)$ qui est nul. Comme le degr\'e de ${\check
K}_{X}(-1)$ est nul et qu'il n'a pas de sections, ce faisceau ne peut
\^etre somme directe de deux faisceaux de degr\'es non nul. Le
faisceau ${\check K}_{X}$ est donc de l'une des deux formes
souhait\'ees (\cite{Ha1} th\'eor\`eme V.2.15).

\vs 0.4 cm

\noi
\th{Proposition}{9} \textit{La vari\'et\'e des quotients inversibles
de degr\'e $3$ de ${\check K}_X$ s'identifie \`a celle des triplets de
points associ\'es \`a $C_0$ sur la cubique $X$.}

\vskip 0.2 cm 

\dm :
La vari\'et\'e des triplets de points associ\'es \`a $C_0$ s'identifie
\`a celle des sections de ${\check K}_{\V}$ :  les sections s'annulent
exactement sur les triplets 
(voir par exemple \cite{Ba1} ou
\cite{P2}). 
Si un triplet $Z$ est sur $X$ alors la
restriction \`a $X$ nous donne une surjection ${\check
K}_{X}\fl\oo_X(2-Z)$, ce dernier faisceau est inversible de degr\'e
$3$. 

R\'eciproquement, soit ${\cal L}$ un quotient inversible de degr\'e 3
de ${\check K}_X$. Notons
$\L'$ le noyau de ${\check K}_X\fl {\cal L}$
qui est inversible de degr\'e $3$. Soit $N$ le
noyau de la compos\'ee ${\check K}_{\V}\fl{\check K}_X\fl\L$. Nous
avons la suite exacte 
$$0\fl{\check K}_{\V}(-3)\fl N\fl\L'\fl 0$$
L'espace vectoriel $H^0N$ est le noyau de la fl\`eche de $H^0\L'$
dans $H^1{\check K}_{\V}(-3)$. Ces groupes sont de dimension
respectives $3$ et $2$. Le faisceau $N$ a donc au moins une section qui nous
d\'efinit une section de ${\check K}_{\V}$ (donc un triplet $Z$) et
ainsi une surjection $\I_Z(2)\fl\L$. La restriction de cette
surjection \`a $X$ nous donne alors le diagramme suivant :
$$\begin{array}{ccccccccccc}
  0 & \fl & \oo_{Z\cap X} & \fl & \I_Z(2)_X & \fl & \oo_X(2-(Z\cap X)) &
    \fl & 0 \\
    &    & &   & \downarrow &   &  &   & \\
   &  &  & &\L &  & &  & \\
\end{array}$$
Comme $\L$ est sans torsion la fl\`eche $\oo_{Z\cap X}\fl\L$ est nulle.
Nous avons donc une fl\`eche surjective de $\oo_X(2-(Z\cap X))$ vers $\L$
qui doit \^etre un isomorphisme. Le degr\'e de $\oo_X(2-(Z\cap X))$
est $6-{\rm{Card}}(Z\cap X)$ et doit \^etre \'egal \`a ${\rm{deg}}(\L)$
c'est \`a dire $3$. Ceci n'est possible que si $Z$ est contenu dans
$X$ c'est \`a dire si la section correspond \`a un triplet sur la
cubique.

Il nous reste \`a v\'erifier que cette section est uniquement
d\'etermin\'ee. Si $h^0N\geq 2$ alors nous avons le diagramme suivant :
$$\begin{array}{ccccccccccc}
  0 & \fl & \oo_{\V}^2 & \fl & {\check K}_{\V} & \fl & \oo_{C_1}(2) &
    \fl & 0 \\
    & & \downarrow & & \downarrow & & \downarrow & & \\
  0 & \fl & \L' & \fl & {\check K}_{X} & \fl & \L & \fl & 0 \\
\end{array}$$
o\`u $C_1$ est une conique. Ceci est absurde car $\L$ est localement
libre sur la cubique $X$ alors que $\oo_{C_1}(2)\fl\L$ doit \^etre
surjective.

\vs 0.4 cm

\th{Remarque}{6}
Cette proposition nous permet de mettre en rapport la proposition
8 et les r\'esultats de \cite{P2}. En effet, nous avons trois cas
selon que ${\check K}_{X}$ a deux, un seul ou une infinit\'e de quotients
inversibles de degr\'e $3$. De la m\^eme fa\c con, nous avons sur $X$
deux, un unique ou une infinit\'e de triplets associ\'es \`a la
conique canonique. Chacun des trois cas se correspondent. 

Parall\`element, les r\'esultats de [P2] nous permettent de dire que
g\'en\'eriquement une cubique a deux triplets en relation avec la
conique canonique, qu'il existe une hypersurface irr\'eductible de
degr\'e $6$ de l'espace des cubiques o\`u il y a un unique triplet et
trois ferm\'es 
irr\'eductibles de codimension $3$ et de degr\'e respectifs $5$, $30$
et $12$ de cubiques ayant une infinit\'e de triplets. Les courbes
g\'en\'erales de ces ferm\'es sont respectivement une cubique
irr\'eductible, la r\'eunion d'une droite et d'une conique en relation
de Poncelet avec $C_0$ et la r\'eunion d'une conique et d'une droite
tangente \`a $C_0$. 


La donn\'ee d'un quotient de rang 3 de 
${\check K}_{X}$
correspond \`a la donn\'ee d'une section de la surface
${\check S}$. C'est une courbe elliptique de degr\'e $3$ c'est \`a
dire une cubique plane. Il y a donc 
selon les cas deux, une unique
ou une infinit\'e de cubiques planes trac\'ees sur ${\check
S}$. L'intersection r\'esiduelle de ${\check S}$ et du plan d'une
telle cubique est form\'ee de trois g\'en\'eratrices de la surface
(les trois points du triplet sur $X$). 
De telles cubiques correspondent donc \`a des points triples
de $S$, il y en a donc deux, un unique ou une infinit\'e.

\vs 0.4 cm

Nous \'etudions maintenant le faisceau $E''$ dans chacun de ces
cas. Remarquons qu'un triplet de points associ\'e \`a $C_0$ ou encore
une section de ${\check K}_{\V}$ correspond \`a un plan $H$ de $\pd$ :
les trois points de $\V$ sont les trois bis\'ecantes \`a $C$ passant
par les points de $H\cap C$. Les triplets correspondent donc aux
points de $\p^3$. L'incidence point/droite restreinte \`a $\V$ est de
degr\'e $3$ au dessus de $\p^3$. C'est l'incidence $I$ entre diviseurs
de degr\'es $2$ et $3$ sur $\pu$ (nous noterons toujours $p$ et $q$ les
projections).

\vs 0.4 cm

\noi
\th{Proposition}{10} \textit{Un point de $\p^3$ est singulier pour $E''$ si et seulement si son triplet associ\'e \`a $C_0$ est sur la cubique $X$. Ces points sont aussi le lieu triple de $S$.}

\vskip 0.3 cm

\dm :
Le faisceau $F$ est un fibr\'e de Scharzenberger pour l'incidence entre $\p^2=\V$ et $\p^3$. La vari\'et\'e d'incidence $I$ a la r\'esolution suivante dans $\p^2\times\p^3$ :
$$0\fl\oo_{\p^2\times\p^3}(-2,-2)\fl K_{\V}\ot\oo_{\p^2\times\p^3}(0,-1)\fl\oo_{\p^2\times\p^3}\fl\oo_I\fl 0$$
ce qui nous donne la r\'esolution :
$$0\fl S_2\ot\oo_{\p^3}(-2)\fl S_1\ot S_4\ot\op(-1)\fl S^2S_3\ot\op\fl p_*q^*\oo_{\p^2}(3)\fl 0$$
Cependant $F(1)$ est le noyau du morphisme de $S_1\ot\op(2)$ dans
$S_4\ot\op(3)$ (cf. proposition 5) donc le complexe d'Eagon Northcott
(voir par exemple [GP]) nous donne la suite exacte suivante (on note
$\oo$ pour $\oo_{\p^3}$) :
$$0\fl S_2\ot\oo(-2)\fl S_1\ot S_4\ot\oo(-1)\fl S^2S_3\ot\oo\fl S_1\ot\oo(2)\stackrel{u}{\fl} S_4\ot\oo(3)\fl 0$$
o\`u $F(1)$ est le noyau de $u$. Le faisceau $p_*q^*\oo_{\p^2}(3)$
s'identifie donc \`a $F(1)$ qui est alors un
module inversible sur la $\oo_{\p^3}$-alg\`ebre finie $p_*\oo_I$. 
La fl\`eche de $\oo_{\p^3}$ dans $F(1)$ est donc nulle si et seulement
si 
celle de $p_*\oo_I$ dans $F(1)$ est
nulle. Mais alors la section de $\oo_{\p^2}(3)$ d\'efinissant $X$ nous
donne le morphisme de $p_*q^*\oo_{\p^2}$ dans $F(1)$ dont le conoyau
est $p_*q^*\oo_X$ (le support de ce faisceau est la surface $S$). Le
lieu triple de ce dernier faisceau d\'ecrit les triplets de points
associ\'es \`a $C_0$ qui sont sur $X$. Ce lieu est donc donn\'e par
l'annulation de la fl\`eche de $p_*\oo_I$ dans $F(1)$ ce qui est
\'equivalent \`a l'annulation de $s$. C'est donc le lieu singulier de
$E''$.

\vskip 0.4 cm

Nous d\'ecrivons dans la proposition suivante le lieu singulier du
faisceau $E''$ lorsqu'il est infini. Ceci a lieu sur trois ferm\'es
irr\'eductibles des cubiques planes (cf. remarque 6).

\vs 0.2 cm

\noi
\th{Proposition}{11} \textit{Si $X$ est une cubique irr\'eductible
ayant une infinit\'e de triplets associ\'es \`a $C_0$, la r\'eunion
d'une droite et d'une conique en relation de Poncelet avec $C_0$
ou la r\'eunion d'une conique et d'une droite tangente \`a
$C_0$ alors le lieu singulier de $E''$ est une cubique gauche, une
droite ou une conique.}

\vskip 0.2 cm

\dm : Les points singuliers de $E''$ sont donn\'es par
les quotients inversibles de rang $3$ de ${\check K}_X$. Or on sait
que dans ce cas ${\check K}_X=\L\oplus\L$ o\`u $\L$ est inversible de
degr\'e $3$. Ainsi les seuls quotients possibles sont isomorphes \`a
$\L$ et ces quotients sont donn\'es par
$\p({\rm{Hom}}(\L\oplus\L,\L))=\pu$. Le lieu singulier de $E''$ est
donc une courbe rationnelle de $\p^3$. Son degr\'e est donn\'e par le
nombre de points dans un hyperplan, c'est \`a dire par le nombre de
triplets tangents en un point fix\'e de la conique. Le degr\'e est 3,
1 et 2 dans chacun des trois cas consid\'er\'es.

\vskip 0.4 cm

Nous r\'esumons
les r\'esultats obtenus sur la famille ${\bf U}$
dans le th\'eor\`eme suivant. Les r\'esultats sur les sous
vari\'et\'es de $\U$ : codimension, irr\'eductibilit\'e et degr\'es
viennent de l'\'etude de la position d'une cubique plane par rapport
\`a une conique fix\'ee d\'etaill\'ee dans 
\cite{P2}.
Notons $\Psi$ le morphisme de $\U$ dans $\H_{3,0}$ qui \`a un faisceau
associe sa courbe des plans instables.

\vskip 0.4 cm

\noi
\th{Th\'eor\`eme}{1} \textit{(\i) Le morphisme $\Psi$ est surjectif. La
fibre de $\Psi$ au dessus de $C\in\H_{3,0}$ est un ouvert de
$\p(H^0\oo_{S^2C}(3))$ la famille des cubiques du plan $S^2C$ (dont
l'image dans $\G$ est la famille des sextiques elliptiques de la
V\'eronese $\V$ des bis\'ecantes \`a $C$).}

\textit{(\i\i) La sous-vari\'et\'e ${\bf U}'$ des faisceaux
r\'eflexifs  ayant un unique point singulier est une hypersurface
irr\'eductible donn\'ee dans chaque fibre par une hypersurface
irr\'eductible de degr\'e 6.}

\textit{(\i\i\i) Le lieu des faisceaux non r\'eflexifs est r\'eunion
de trois ferm\'es irr\'eductibles de codimension 3
donn\'es dans les fibres par des ferm\'es irr\'eductibles de
degr\'es $5$, $30$ et $12$.}

\vs 0.4 cm

Notons $\H^d_{6,3}$ (resp. $\H^2_{6,3}$) le localement ferm\'e de
$\H_{6,3}$ (resp. $\H^d_{6,3}$) des courbes r\'eunion d'une cubique
gauche irr\'eductible et d'une cubique plane qui se coupent en trois
points (resp. dont la surface des tris\'ecantes non contenues dans le
plan de la cubique porte une seule cubique plane). 

Nous pouvons d\'efinir une application rationnelle $r$ de $\H^d_{6,3}$
vers $\U$ : la cubique gauche $C$ d\'efinit un point de $\H_{3,0}$ et la
cubique plane d\'etermine une courbe elliptique de degr\'e $6$ sur la
V\'eronese $\V$ des bis\'ecantes \`a $C$. En effet, \`a un point
de la cubique plane on associe la bis\'ecante \`a la cubique gauche
passant par ce point. Le th\'eor\`eme 1 nous permet de dire que ceci
d\'efinit un \'el\'ement $E''$ de $\U$. La cubique gauche $C$ est
alors la courbe des plans instables de $E''$.


\vs 0.4 cm

\noi
\th{Corollaire}{3} \textit{L'application rationnelle $r$ est
dominante de degr\'e $2$. Elle est ramifi\'e au dessus de
$\U'$. L'image r\'eciproque de $\U'$ est $\H^2_{6,3}$.}

\vs 0.2 cm

\dm :
Une courbe de la fibre au dessus de $E''$ est
r\'eunion d'une cubique gauche $C$ et d'une cubique plane $C'$ la
rencontrant en trois points. La cubique gauche $C$ est la courbe des
plans instables de $E''$. La surface des
trisc\'ecantes \`a $C\cup C'$ (non contenues dans le plan de $C'$) est la
surface $S$ recouverte par les
droites bisauteuses de $E''$. La courbe $C$ et la surface $S$ sont
donc fix\'ees par $E''$. La fibre est d\'etermin\'ee par les cubiques
planes $C'$ trac\'ees sur $S$.
%
%
L'\'etude pr\'ed\'edente nous permet de dire que g\'en\'eriquement il
y a deux telles cubiques et qu'au dessus de $\U'$ il y en a une
unique. La fibre au dessus des trois composantes du lieu o\`u le
faisceau est non r\'eflexif est infinie.

\subsection*{Lien avec les instantons de degr\'e 3}

Nous montrons le lien entre la famille $\U$ et 
les instantons de degr\'e 3. Donnons nous, pour un faisceau $E''\in\U'$ dont
le lieu singulier est concentr\'e au point $P$, une
fl\`eche surjective vers le faisceau du point : $E'' \fl {\cal
O}_P$. Le noyau $E$ est un \'el\'ement de $ M_{\p^3}(0,3,0)$. Par
l'interm\'ediaire de la vari\'et\'e $\F$, nous d\'eterminons 
un morphisme canonique $E''\fl\oo_P$ dont le noyau est dans le bord de
$\bI_3$. Ces faisceaux forment une composante irr\'eductible du bord. 

Pour construire les d\'eformations nous 
utilisons les transformations
cubo-cubiques. Nous avons vu que 
que la famille $\di$ et la famille $\U'$ sont
chacune birationnelles \`a des ferm\'es du sch\'ema de $\H_{6,3}$. Le
morphisme $f$ nous permet de remonter ces courbes en des
transformations cubo-cubiques. Notons $\T_3^t$, $\T_3^1$, $\T_3^d$
et $\T_3^2$ les images r\'eciproques par $f$ de $\H_{6,3}^t$,
$\H_{6,3}^1$, $\H_{6,3}^d$ et $\H_{6,3}^2$.

\vs 0.4 cm

\noi
\th{Proposition}{12} \textit{Les transformations cubo-cubiques de
$\T_3^d$ (resp. $\T_3^2$) sont obtenues comme inverses de
celles de $\T_3^t$ (resp. $\T_3^1$). La vari\'et\'e
$\T_3^1$ (resp. $\T_3^2$) est l'intersection de $\pi(\F)$ avec
$\T_3^t$ (resp. $\T_3^d$). Enfin on a $\pi^{-1}(\T_3^1)=\F_3$ et
$\pi^{-1}(\T_3^2)=\F_1$.}

\vs 0.2 cm

\dm :
Nous commen\c cons par d\'ecrire les r\'esolutions des id\'eaux des
courbes des diff\'erentes sous-vari\'et\'es du sch\'ema de Hilbert.

\vs 0.2 cm

\noi
\th{Lemme}{3} \textit{Soit $t\in\T_{3}^t$, il existe des bases de $W$
et de $V$ telles que cette application s'\'ecrive comme une matrice
$4\times 4$ (de $V$ dans $W$) \`a coefficients dans ${\check R}$ sous
la forme suivante :
$$\left(\begin{matrix}
 A & 0 \\
 C & D \\
\end{matrix}\right)$$
o\`u $A$ est de taille $3\times 3$. De plus, on peut choisir $A$
sym\'etrique si et seulement si $t\in\T_3^1$.}

\vskip 0.2 cm

\dm :
Soit $t\in\T_3^t$ et soit $Y$ la courbe de $\H_{6,3}^t$ correspondant
\`a $t$. Son id\'eal est r\'esolu par :
$$0 \fl R\otimes \oo_{\pvv}(-4)\stackrel{t}{\fl} W\otimes\oo_{\pvv}(-3)\fl \I_Y\fl 0$$
Le point triple est donn\'e par le deuxi\`eme id\'eal de Fitting de
$\I_Y$ donc par les mineurs $2 \times 2$ de la matrice $3\times 4$ (de
$R$ dans $W$) \`a coefficients dans ${\check V}$. Ils sont tous dans
l'id\'eal engendr\'e par le noyau de $H\in V$ d\'efinissant le
point triple. 
Nous avons donc une matrice $4\times 4$ (de $V$ dans $W$) \`a coefficients
dans ${\check R}$ de la forme voulue. Si une
transformation cubo-cubiques peut s'\'ecrire sous cette forme elle
donne \'evidement une courbe de degr\'e $6$ et de genre $3$ ayant un
point triple.

De plus, la courbe $Y$ est sur un complexe de droites si et seulement si la
matrice est la multiplication du modules de Rao d'un faisceau de $\di$
ce qui nous donne exactement la condition de sym\'etrie (voir le
paragraphe pr\'ec\'edent : \'etude de $\di$).

\vskip 0.2 cm

\noi
\th{Lemme}{4} \textit{Soit $t\in\T_3^d$, il existe des bases de $W$ et de $V$ telles que cette application
s'\'ecrive comme une matrice $4\times 4$ (de $V$ dans
$W$) \`a coefficients dans ${\check R}$ sous la forme
suivante :
$$\left(\begin{matrix}
 A & 0 \\
 C & D \\
\end{matrix}\right)$$
o\`u $D$ est de taille $3\times 3$. De plus, on peut choisir $D$
sym\'etrique si et seulement si $t\in\T_3^2$.}

\vskip 0.2 cm 

\dm :
Soit $t\in\T_3^d$ et soit $Y'$ la courbe correpondante qui est ACM de
degr\'e $6$ et de genre $3$. Son id\'eal v\'erifie :
$$0\fl \I_{Y'} \fl \I_C \fl {\cal L}' \fl 0$$
o\`u $C$ est une cubique gauche et ${\cal L}'$ est support\'e par une
cubique plane $C'$ rencontrant $C$ en trois points. Les r\'esolutions
de $\L'$ et $\I_C$ sont :
$$0\fl R\otimes \oo_{\pvv}(-4) \fl (R\oplus B)\otimes\oo_{\pvv}(-3)
\stackrel{(a,H)}{\fl} B\otimes \oo_{\pvv}(-2) \fl {\cal L}' \fl 0$$
et
$$0\fl R_0\otimes \oo_{\pvv}(-3) \fl {W_0}\otimes \oo_{\pvv}(-2)
\fl \I_C \fl 0$$
o\`u $R$ et $B$ sont des espaces vectoriels de dimension 3 et $R_0$ et
${W_0}$ sont des espaces vectoriels de dimensions respectives 2 et
3. La fl\`eche de $\I_C$ dans $\L'$ nous d\'efinit des fl\`eches de
$W_0$ dans $B$ et de $R_0$ dans $R\oplus B$. Comme $Y'$ est ACM, on a
n\'ecessairement $W_0\backsimeq B$. Ceci impose que l'application de $R_0$ dans
$R\oplus B$ est injective. Notons $W$ le conoyau de cette injection,
nous avons la r\'esolution de l'id\'eal $\I_{Y'}$ :
$$0\fl R\otimes \oo_{\pvv}(-4) \stackrel{M}{\fl} W\otimes
\oo_{\pvv}(-3) \fl \I_{Y'}\fl 0$$
Enfin, comme $C$ n'a pas de composante dans le plan de $C'$, on voit,
par restriction \`a ce plan, que $R_0\cap B$ est nul dans $R\oplus
B$. Ceci nous permet de dire que les fl\`eches $R_0\fl R$ et $W_0\fl
W$ sont injectives. En prenant des bases de $R_0$ et $W_0$
compl\'et\'ees en des bases de $R$ et de ${W}$, la matrice $M$
s'\'ecrit :
$$\left(\begin{matrix}
 N & * \\
 0 & H \\
\end{matrix}\right)$$
o\`u $N$ est la matrice de la cubique gauche et $H$ est l'\'equation
du plan de $C'$. Soit $V_0$ le sous-espace vectoriel de dimension $3$
de $V$ qui est le noyau de la forme lin\'eaire $H\in\cv$. En prenant
un base de $V_0$ compl\'et\'ee en une base de $V$ on a une matrice de
la forme voulue. Dans les notations de l'\'enonc\'e, la matrice $D$
repr\'esente 
l'\'el\'ement de ${\check R}\otimes {\check V_0}\otimes {W_0}$. Si une
transformation cubo-cubique s'\'ecrit sous cette forme elle donne une
courbe de $\H_{6,3}^d$.

La courbe est dans $\H_{6,3}^2$ si il y a une unique cubique plane sur
la surface des tris\'ecantes. Or nous avons vu que ceci n'est possible que
si le faisceau $K(1)_X$ est extension non triviale d'un faisceau
${\cal L}$ par lui m\^eme qui est alors une th\'eta-caract\'eristique
(car $\Lambda^2 K(1)_X=\oo_X$). La courbe $X$ est la courbe des
droites param\'etrisant les tris\'ecantes vue dans le plan 
des bis\'ecantes de $C$. 

De plus, l'application
rationnelle du plan de $X$ dans celui de $C'$ qui a une bis\'ecante de
$C$ associe son point d'intersection avec le plan de $C'$ est la
transformation quadratique associ\'ee aux sommets du triangle de
Poncelet sur $X$ (c'est \`a dire aux points d\'efinis par $\L$). Sa
r\'eciproque est la transformation quadratique associ\'ee aux points
de $C\cap C'$. Cette application induit une bijection de $X$ sur
$C'$. Si ${\cal L}'$ est donn\'e par la
multiplication $R\otimes\oo_{{\p}({\check V_0})}(-3)\fl
{W_0}\otimes \oo_{{\p}({\check V_0})}(-2)$,
alors ${\cal L}$ est donn\'e par la
multiplication $R\otimes \oo_{{\p}({W_0})}(-3)\fl {\check
V_0}\otimes \oo_{{\p}({W_0})}(-2)$. Le faisceau ${\cal L}$ \'etant
une th\'eta-caract\'eristique, il existe un isomorphisme entre
${\check R}$ et ${\check V_0}$ tel que la compos\'ee $V_0\ot V_0\fl
R\ot V_0\fl W_0$ se factorise par $S^2V_0$.

Nous utilisons alors le r\'esultat de Barth \cite{Ba1} suivant : si on a trois
espaces vectoriels de dimension trois  $R_0$, ${\check V_0}$ et $W_0$ et
un \'el\'ement de ${\check R}\otimes {\check V_0}\otimes { W_0}$ tel qu'il existe
un isomorphisme entre ${\check R}$ et ${\check V_0}$ tel que la compos\'ee $V_0\ot V_0\fl R\ot V_0\fl W_0$ se factorise par $S^2V_0$, alors il existe deux autres isomorphismes entre $W_0$ et ${\check V_0}$ et entre ${\check R}$ et
${W_0}$ tels que l'application $R\fl {\check V_0}\otimes {\check V_0}$
(resp. ${\check V_0} \fl {W_0}\otimes {W_0}$) se factorise par
$S^2{\check V_0}$ (resp. $S^2W_0$). Ceci n'est vrai qu'en
dimension $3$ car alors tous les r\'eseaux de quadriques v\'erifient
la condition $\a_3$ de Barth (cf. \cite{Ba1}). Ce r\'esultat nous
permet de choisir $A$ sym\'etrique si et seulement si
$t\in\T_3^2$.

\vs 0.2 cm

Les transformations
cubo-cubiques de $\T_3^d$ (resp. $\T_3^2$) sont les inverses de celles
de $\T_3^t$ (resp. $\T_3^1$) car elles sont obtenues par \'echange des
r\^oles de $\cv$ et $W$. De plus, on sait que $\T^1_3$ est l'image de
$\F_3$ dans $\T^t_3$ (proposition 2) ainsi, par \'echange des
r\^oles, nous voyons que $\T^2_3$ est l'image de $\F_1$ dans
$\T^d_3$. Nous avons ici disym\'etris\'e les r\^oles de $W$ et $\cv$
ce qui nous donne la \textit{dualit\'e} entre ces deux cas.

\vs 0.4 cm

\noi
\th{Proposition}{13} \textit{Le complexe d\'ecrit en introduction
d\'efinit un morphisme $h$ d'un ouvert de $\T_3^d$ vers $\Ma$. Ce
morphisme est la compos\'ee de $f$ avec $r$. L'image r\'eciproque de
$\U'$ est $\T_3^2$.}

\vs 0.2 cm

\dm :
Nous construisons en fait un morphisme vers $\Ma$ \`a
partir d'une vari\'et\'e dominant birationnellement $\T_3^d$ ce qui
nous donnera notre morphisme sur un ouvert. Consid\'erons ainsi la
sous-vari\'et\'e $\F^1$ 
de $\T\times\p({\rm{Hom}}(W,{\check V}))$.
Cette vari\'et\'e contient la projection naturelle de $\F_1$ comme une
sous-vari\'et\'e de codimension $1$. De plus $\F^1$ domine
birationnellement $\T_3^d$. En effet, le morphisme de $\F^1$ vers $\T$
est birationnel sur son image (fait 3) qui est le  $\T_3^d$ (lemme
$4$). 
Sur la vari\'et\'e $\F^1$ nous
avons un complexe universel :
$$R\otimes \Omega^2(2) \stackrel{\varphi}{\fl} W\otimes \Omega^1(1)
\stackrel{\psi}{\fl} \oo_{{\p}(V)}$$
qui est injectif sur un ouvert de $\F^1$ et jamais surjectif. Sa
cohomologie au centre nous donne un faisceau plat qui d\'ecrit un
morphisme vers $\Ma$.

Soit $t\in\T_3^d$, son image par $f$ nous donne une courbe $Y'$ de
$\H^d_{6,3}$. Le lemme $4$ nous permet d'\'ecrire le diagramme suivant :
$$\begin{array}{ccccccccccc}
& & & & 0 & & 0 & &  & & \\
& & & & \downarrow & & \downarrow & & & & \\
 & & 0 & \fl & R_0\otimes \om^2(2) & \fl  & W_0\otimes\om^1(1) & \fl & F &
 \fl & 0 \\
 & & & & \downarrow & & \downarrow & & & & \\
& & 0 & \fl & R\otimes \om^2(2) & \fl & W\otimes\om^1(1) & \fl & Q &
 \fl & 0 \\
 & & & & \downarrow & & \downarrow & & & & \\
 0 & \fl & \oo_{\p^3}(-1) & \fl & \om^2(2) & \stackrel{v}{\fl} & \om^1(1) & \fl & \I_P
 & \fl & 0 \\
& & & & \downarrow & & \downarrow & & & & \\
& & & & 0 & & 0 & & & & \\
\end{array}$$
o\`u le faisceau $F$ est le fibr\'e de Schwarzenberger associ\'e \`a
la cubique gauche (cf. proposition 5) et $v$ est l'\'equation du
plan $H$ de $\pd$ qui correspond \`a un point $P$ de $\p^3$. Le lemme
du serpent nous donne la suite exacte suivante $0 \fl \oo_{\p^3}(-1)
\fl F \fl Q \fl \I_P \fl 0$. La construction par le complexe nous dit
que $E''(1)$ est le noyau de $Q\fl\I_P$ alors que nous savons que
l'image de $r$ est donn\'ee par le conoyau de $\oo_{\p^3}(-1)\fl
F$.

\vs 0.4 cm

Nous calculons maintenant la limite d'une famille d'instantons dont le
terme dans $\F$ tend vers un \'el\'ement $t_0$ de $\F_1$. 
Soit donc $A$ un anneau de
valuation discr\`ete de corps r\'esiduel $k$ et $t_A$ un $A$-point de
$\F$ dont le point g\'en\'erique $t_g$ est dans $\F_4$ et le point
sp\'ecial $t_0$ est assez g\'en\'eral dans $\F_1$. Soit ${\cal E}$ la
famille plate obtenue \`a partir de $t_A$ et du morphisme $g$.

\vs 0.4 cm

\noi
\th{Proposition 14} \textit{Le bidual $E''$ du faisceau limite
$E$ de ${\cal E}$ est le faisceau de $\Ma$ donn\'e par
$h\circ\pi(t_0)$.}


\vs 0.2 cm

\dm : 
On construit la monade $R\ot\om^2(2)_A\fl W\ot\om^1(1)_A\fl\I_{A.P}$
(o\`u $\I_{A,P}$ est l'id\'eal dans $\p^3_A$ du point $P\in\p^3_k$ )
obtenue \`a partir du point $t_A$. Le faisceau ${\cal E}$ est la
cohomologie au centre de cette monade. Les hypoth\`eses de
g\'en\'ericit\'e nous permettent de dire que la premi\`ere fl\`eche
horizontale est injective m\^eme au point sp\'ecial ce qui nous permet
de conclure \`a la platitude de ${\cal E}$. La monade au point
sp\'ecial est :
$$\begin{array}{ccccccccc}
   &   &   &   & 0 &   & 0 &   &  \\
   &   &   &   & \downarrow &   & \downarrow &   &  \\
 0  & \fl  & R \otimes \Omega^2(2) & \fl &  K & \fl &
    E(1) & \fl & 0 \\
   &   & \Vert &   & \downarrow &   & \downarrow &   &  \\
 0  & \fl  & R \otimes \Omega^2(2) & \fl & W \otimes
   \Omega^1(1) & \fl &  Q & \fl & 0 \\
   &   &   &   & \downarrow &   & \downarrow &   &  \\
   &   &   &   & \L & = & \L &   &  \\ 
   &   &   &   & \downarrow &   & \downarrow &   &  \\
   &   &   &   & 0 &   & 0 &   &  \\
\end{array}$$
o\`u $\L$ est une extension de $\I_P$ par $\oo_P$. Le bidual de $E(1)$ est \'evidement donn\'e par le noyau de la fl\`eche $Q\fl\I_P$ d\'eduite de $Q\fl\L$ ce qui nous dit qu'il est donn\'e par $h\circ\pi(t_0)$.

Remarquons que le faisceau $\L$ est une extension triviale
c'est \`a dire est isomorphe \`a $\I_P\oplus\oo_P$. En effet, cette
extension est donn\'ee par $\I_{A,P}\ot k$.

\vs 0.4 cm

Soit $E''\in\U$, soient $\vp$ et $\psi$ au dessus de $E''$ dans
$\F^1$. Ces morphismes correspondent au choix d'un \'el\'ement dans la
fibre de $r$ (il y a, \`a priori, deux choix mais si $E''\in\U'$ il n'y
aura pas d'ambiguit\'e) puis d'un \'el\'ement de $PGL(R)\times
PGL(W)$.
Soit $Q$ le conoyau de la fl\`eche
$R\ot\om^2(2)\stackrel{\vp}{\fl} W\ot\om^1(1)$. Le faisceau $E''$ est
le noyau de la surjection $Q\fl\I_P$ d\'eduite de $\psi$. Nous avons donc
une fl\`eche $\hom(E'',\oo_P)\fl \ext(\I_p,\oo_P)$. La proposition
pr\'ec\'edente nous dit que pour qu'il existe une fl\`eche $s$ de
$E''$ dans $\oo_P$ dont le noyau est un faisceau du bord de $\bI_3$,
il faut que l'image de $s$ dans $\ext(\I_p,\oo_P)$ soit nulle. Notons
$(*)$ cette condition.


\vs 0.4 cm

\noi
\th{Proposition}{15} \textit{Un faisceau $E''\in\U$ v\'erifie la
condition $(*)$ si et seulement si $E''\in\overline{\U'}$. Sur $\U'$,
il existe un unique \'el\'ement de $\hom(E'',\oo_P)$ nul dans
$\ext(\I_p,\oo_P)$. Cet \'el\'ement nous d\'efinit un morphisme $\a$
de $\U'$ dans $\M$ birationnel sur son image.}

\vs 0.2 cm

\dm : 
Le faisceau $E''$ est le noyau de la surjection $Q\fl\I_P$ d\'eduite
de $\psi$. Le point $P$ d\'etermine un sous-espace vectoriel $V_0$ de
dimension $3$ de $V$ et $\psi$ d\'etermine un sous-espace vectoriel
$W_0={\rm{Ker}}\psi$ de dimension $3$ de $W$ et un quotient $W_1$ de
rang $1$. On a la suite exacte :
$$0\fl {\check
V_0}\fl\hom(Q,\oo_P)\fl\hom(E'',\oo_P)\fl\ext(\I_p,\oo_P)$$
Nous allons montrer que selon que $\psi'$ existe ou non (i.e. selon que
$E''\in\overline{\U'}$ ou non), l'espace $\hom(Q,\oo_P)$ est de
dimension $4$ ou $3$ et qu'ainsi il existe ou non un \'el\'ement de
$\p({\rm{Hom}}(E'',\oo_P))$ qui donne une extension triviale de $\I_p$
par $\oo_P$.

La d\'efinition de $Q$ nous permet de donner la suite exacte suivante :
$$0\fl\hom(Q,\oo_P)\fl{\check W}\ot{\check V_0}\fl{\check R}\ot V_0\fl\ext(Q,\oo_P)\fl 0$$
o\`u la fl\`eche centrale est la compos\'ee :
$${\check W}\ot{\check V_0}\stackrel{\vp}{\fl}{\check
R}\ot\cv\ot{\check V_0}\fl{\check R}\ot{\check V_0}\ot{\check
V_0}\fl{\check R}\ot\Lambda^2{\check V_0}\fl{\check R}\ot V_0$$
L'espace vectoriel ${\check W_1}\ot{\check V_0}$ est donc contenu dans
$\hom(Q,\oo_P)$. Pour savoir si le noyau est r\'eduit \`a cet espace
ou est de dimension sup\'erieure \`a trois il faut \'etudier la
fl\`eche r\'eduite qui est la compos\'ee : ${\check W_0}\ot{\check
V_0}\fl{\check R}\ot{\check V_0}\ot{\check V_0}\fl{\check
R}\ot\Lambda^2{\check V_0}\fl{\check R}\ot V_0$. Cette fl\`eche est
non bijective si l'image de ${\check W_0}\ot{\check V_0}$ dans
${\check R}\ot{\check V_0}\ot{\check V_0}$ rencontre ${\check R}\ot
S^2{\check V_0}$ c'est \`a dire exactement si on peut trouver un
morphisme de $W_0$ dans ${\check V_0}$ qui rend la fl\`eche $\vp$
sym\'etrique. Ceci est exactement la condition d'appartenance \`a
$\overline{\U'}$. Si le faisceau $E$ est dans l'image de $\alpha$,
alors on retrouve $E''$ qui est le bidual de $E$.

Remarquons qui si on est dans un des trois ferm\'es de codimension $3$
de $\U$ pour lesquels $E''$ est non r\'eflexif, alors le choix du
morphisme de $E''$ dans $\oo_P$ n'est plus canonique (on a un $\pu$
d'\'el\'ements qui donnent une extension triviale de $\I_P$ par
$\oo_P$). Le morphisme $\alpha$ n'est donc pas d\'efini sur ces
ferm\'es.

\vs 0.4 cm


Nous montrons enfin que le morphisme obtenu par composition de
$h\circ\pi$ de $\F_1$ dans $\U'$ avec $\alpha$ d\'efinit un
morphisme de $\F_1$ dans $\M$ qui prolonge $g$ sur $\F_1$. Pour ceci,
nous montrons que pour toute famille d'instantons ${\cal E}$ construite
comme pour la proposition 14, la limite $E$ est l'image de $t_0$
par $\alpha\circ h\circ\pi$.

%
%
%
%

\vs 0.4 cm

\noi
\th{Corollaire}{4} \textit{L'image de $\F_1$ par $\alpha\circ
h\circ\pi$ est une vari\'et\'e irr\'eductible $\dj$ de $\M$ qui est
dans le bord de $\bI_3$. Ce morphisme prolonge le morphisme $g$ sur
$\F_1$.}

\vs 0.2 cm

\dm :
Soit ${\cal E}$ une famille d'instantons construite comme pour la
proposition 14. Soit $E$ le faisceau limite et $E''$ son
bidual. Nous avons vu que $E''$ est $h\circ\pi(t_0)$ et
que $E$ est le noyau d'une fl\`eche $s\in\hom(E'',\oo_P)$ qui s'annule
dans $\ext(\I_P,\oo_P)$. Mais alors ceci impose (proposition 15)
que $E''\in\U'$ (car $t_0$ est assez g\'en\'eral) et que le faisceau
$E$ est exactement $\a(E'')$. 

R\'eciproquement, si $E$ est g\'en\'eral dans l'image de $\a$, alors
son bidual $E''$ est g\'en\'eral dans $\U'$ et l'\'el\'ement $t_0$
correspondant peut \^etre choisit g\'en\'eral. Nous pouvons alors
construire la d\'eformation comme pour la proposition 14 dont la
limite est $E$.

Remarquons que ce morphisme est d\'efini par un faisceau universel sur
$\F_1$. En effet, sur $\F_1$, le morphisme $h\circ\pi$ \`a valeur dans
$\Ma$ est d\'efini par le faisceau universel ${\bf E''}$ donn\'e sur
$\F_1\times\p^3$ par la cohomologie au centre du complexe
$$R\ot\om^2(2)\fl W\ot\om^1(1)\fl\oo$$
Au dessus du localement ferm\'e
$\U'$ de $\U$, si on note $Z$ le ferm\'e $\U'\times\p^3$ d\'efini par
le lieu singulier, le noyau de la fl\`eche de fibr\'es
${\cal H}om({\bf E''},\oo_Z)\fl{\cal E}xt^1(\I_Z,\oo_Z)$ est
inversible. Ce faisceau nous donne une section de $\U'$ dans le
fibr\'e ${\cal H}om({\bf E''},\oo_Z)$ et le noyau ${\bf E}$ de la
fl\`eche universelle ${\bf E''}\fl\oo_Z$ est le faisceau recherch\'e.

\vs 0.4 cm

\th{Remarque}{7} Nous d\'ecrivons ici plus pr\'ecis\'ement la
conjecture de L. Gruson et G. Trautmann. Nous donnons \'egalement
quelques r\'esultats compl\'ementaires sur cette conjecture.

La 
conjecture affirme que tous les faisceaux $E$ du bord de $\bI_3$ sont
sans torsion non localement libres. Deux cas se pr\'esentent selon que
le bidual $E''$ du faisceau $E$ est localement libre ou non. Dans le
premier cas (type $1$), le lieu singulier de $E$ est une courbe de
degr\'e inf\'erieur ou \'egal \`a $3$. Dans le second cas (type $2$),
le bidual est seulement r\'eflexif et on conjecture que le lieu
singulier d'un faisceau g\'en\'eral $E$ est concentr\'e en un
point. Le bord de $\bI_3$ aurait alors quatre composantes de chaque
type. 

Les courbes associ\'ees aux faisceaux de type $1$ seraient : une
droite, une cubique gauche (ces deux cas apparaissent dans \cite{GS}), une
conique (c'est le cas de $\di$ que nous traitons ici) et une cubique
plane (ce cas peut \^etre obtenu par ``d\'eformation par
homoth\'etie'').
%
%

Les quatres composantes de type $2$ correspondraient \`a des faisceaux
r\'eflexifs de troisi\`eme classe de Chern vallant $2$, $4$, $6$ ou
$8$. Nous traitons ici le cas $c_3=2$ (vari\'et\'e $\dj$). Des
faisceaux appartenant aux cas $c_3=4$ ou $6$ peuvent \^etre obtenus
par ``d\'eformation de courbes''
%
et pour $c_3=8$ par ``d\'eformation par homoth\'etie''.
%
%

Concernant la conjecture de L. Gruson et G. Trautman
sur le bord de $\bI_3$, nous savons construire sept composantes
irr\'eductibles en codimension $1$ sur les huit pr\'edites. L'existence
de la composante du second type correspondant \`a $c_3=6$ est plus
incertaine : nous savons construire des familles de tels faisceaux mais qui
sont en codimension au moins 2.

\section{Quelques situations g\'eom\'etriques associ\'ees}

Nous d\'ecrivons dans ce paragraphe deux applications g\'eom\'etriques
des r\'esultats pr\'ec\'edents. Nous donnons une param\'etrisation de
l'espace des modules des courbes de degr\'e $7$ et de genre $2$. Nous
exhibons une famille $\mathfrak{I}$ de dimension 36 d'involutions
birationnelles de $\p^3$ et nous donnons une repr\'esentation
birationnelle de cette famille. 

\subsection*{Espace des modules des courbes de degr\'e $7$ et de genre
$2$}

Notons $\mathfrak{H}_{7,2}$ le sch\'ema de Hilbert des courbes de
degr\'e $7$ et de genre $2$ de $\p^3$.
\def \hi {\mathfrak{H}_{7,2}}
\def \ho {\mathfrak{H}_{3,0}}
Nous avons vu dans l'exemple $1$ que si $Y$ est une courbe de degr\'e $7$
et de genre $2$ irr\'eductible et lisse alors le groupe
${\rm{Ext}}^1_{{\cal O}_{\p^3}}({\cal I}_Y(2), {\cal
O}_{\p^3}(-2))=H^0\omega_Y$ est de dimension $2$. Ainsi il existe un
pinceau d'extensions non nulles. Chacune de ces extensions donne un
faisceau $E''\in\U_0$ qui est r\'eflexif. La remarque $4$ nous
permet alors de dire que le faisceau $F$ d\'eduit de $E''$ gr\^ace \`a
la suite spectrale de Beilinson est localement libre et donc
associ\'e \`a une cubique gauche irr\'eductible
$C$ (proposition $5$). Le faisceau $E''$ est alors dans l'ouvert
$\U$. La fl\`eche
naturelle de $F(1)$ dans $\I_Y(4)$ (qui se factorise par $E''(2)$)
nous donne la suite exacte :
$$0\fl\oo_{\p^3}^2\fl F(1)\fl\I_Y(4)\fl 0 \ \ \  (*)$$
Cette construction est ind\'ependante du choix de l'extension et nous
permet de d\'efinir un morphisme $\Phi$ de $\hi$ dans $\ho$ le
sch\'ema de Hilbert des cubiques gauches irr\'eductibles qui a la
courbe $Y$ associe la cubique gauche $C$ qui d\'efinit $F$.

\vs 0.4 cm

\th{Proposition}{16} \textit{La fibre de $\Phi$ au dessus d'une
courbe $C\in\ho$ est birationnellement isomorphe \`a
$\G(2,H^0\oo_{S^2C}(3))$ la vari\'et\'e des pinceaux de cubiques du
plan $S^2C$.}

\vs 0.2 cm

\dm :
Soit donc $C\in\ho$. La donn\'ee de $Y\in\hi$ au dessus de $C$ nous
permet de d\'efinir un \'el\'ement de $\G(2,H^0\oo_{S^2C}(3))$ gr\^ace
\`a la suite exacte $(*)$. R\'eciproquement, soit un morphisme
$\oo_{\p^3}^2\fl F(1)$ 
donn\'e par un \'el\'ement de $\G(2,H^0\oo_{S^2C}(3))$. Fixons un point
de ce pinceau $\oo_{\p^3}\fl F(1)$. Le conoyau est alors un faisceau
$E''\in\U$ et on a une suite exacte :
$$0\fl\oo_{\p^3}\fl E''(2)\fl\I_Y(4)\fl 0$$
qui nous d\'efinit une courbe $Y$ de degr\'e $7$ et de genre $2$ qui
est ind\'ependente du point du pinceau fix\'e. Le
faisceau $E''(2)$ est r\'egulier au sens de Castelnuovo-Mumford et
pour un pinceau g\'en\'eral, la courbe $Y$ est lisse.

\vs 0.4 cm

\th{Remarque}{8} L'espace des modules des courbes de degr\'e $7$ et de
genre $2$ de $\p^3$ est le quotient par $PGL_4$ de $\H_{7,2}$. Le
sch\'ema $\ho$ n'a qu'une orbite sous $PGL_4$ (il est isomorphe au quotient $PGL_4/PGL_2$) donc l'espace des modules est birationnellement isomorphe au
quotient par $PGL_2$ de la vari\'et\'e $\G(2,S^3S_2)$ des pinceaux de
cubiques du plan $\p(S_2)$. 

Nous pouvons \`a partir du pinceau retrouver le
mod\`ele non singulier de la courbe : nous avons une droite $\pu$ dans
$\p(S^3S_2)$, la courbe $Y$ est un rev\^etement double au dessus. En
effet, si nous choisissons un point de $\pu$, c'est \`a dire une section de
$F(1)$, elle d\'etermine $E''(2)$ et donc deux points de $Y$
donn\'es par les deux points singuliers de $E''$. Ce sont les deux
points o\`u la section de $\omega_Y$ d\'efinissant $E''$
s'annule. Quand on fait varier le point de $\pu$, on fait varier le
faisceau $E''$ c'est \`a dire l'\'el\'ement de $H^0\omega_Y$. On
recouvre de cette fa\c con tout $Y$. Le rev\^etement double de $Y$ au
dessus de $\pu$ est donn\'e par le ${\bf g}^1_2$ d\'efini par
$\omega_Y$. Les points de ramification de ce morphisme apparaissent
quand la section de $\omega_Y$ s'annule doublement en un point. Ceci
correspond exactement aux faisceaux de $\U'$ ou encore aux cubiques
ayant un unique triplet associ\'e \`a la conique canonique. L'image
des points de ramification est donc form\'ee par l'intersection du
pinceau $\pu$ avec la vari\'et\'e des cubiques ayant un unique
triplet. La formule d'H\" urwitz nous permet de retrouver le fait que
le degr\'e de cette derni\`ere vari\'et\'e est $6$. La courbe
abstraite param\'etrisant $Y$ est le rev\^etement double de $\pu$
ramifi\'e aux $6$ points d'intersection du pinceau et de cette
vari\'et\'e.

\subsection*{Involutions birationnelles de $\p^3$}

Les involutions du plan projectif $\p^2$ sont g\'eom\'etriquement
connues depuis longtemps (voir par exemple [Be]) mais la preuve
rigoureuse de leur classification est plus r\'ecente. A. Beauville et
L. Bayle [BB] ont donn\'e une preuve simplifi\'ee de ce r\'esultat en
utilisant la th\'eorie de Mori. Dans $\p^3$, il ne semble pas qu'il
existe, comme dans le cas de $\p^2$, une classification des
involutions birationnelles. Je d\'ecris ci-dessous une famille de
telles involutions.

\vs 0.4 cm

\th{Fait}{5} \textit{Toute courbe de degr\'e 9 et de genre 6 d\'efinit
une involution birationnelle de $\p^3$. La famille d'involutions
$\mathfrak{I}$ ainsi construite est birationnelle au sch\'ema
$\H_{9,6}$ des courbes de degr\'e 9 et de genre 6.}

\vs 0.2 cm

\dm :
Consid\'erons une courbe $X$ de degr\'e $9$ et de genre $6$
g\'en\'erale, elle est toujours sur quatre quartiques
(i.e. $H^0\I_X(4)$ est de dimension 4). Prenons trois de ces
quartiques, l'intersection r\'esiduelle \`a $X$ est form\'ee de deux
points. Nous d\'ecrivons ainsi une famille de dimension 3
(birationnelle \`a $\p(H^0\I_X(4))$ ce qui d\'efinit une involution.
 
La courbe de degr\'e 9 et de genre 6 est le lieu o\`u l'involution
n'est pas d\'efinie. Le sch\'ema de Hilbert $\H_{9,6}$ des courbes de
degr\'e 9 et de genre 6 param\'etrise donc birationnellement cette
famille d'involutions. Nous voyons  ainsi apparaitre une famille
\`a $36$ param\`etres d'involutions de $\p^3$.

\vs 0.4 cm

\th{Remarques}{9} \textit{(\i)} Ces involutions sont li\'ees \`a notre
situation de la fa\c con suivante : prenons deux quartiques contenant
$X$ (i.e. un sous espace vectoriel $U$ de dimension 2 de
$H^0\I_X(4)$), l'intersection r\'esiduelle est une courbe $Y$ de
degr\'e $7$ et de genre $2$. Lorsqu'on fait varier une troisi\`eme
quartique dans le pinceau $\p(H^0\I_X(4)/U)$, celle-ci d\'ecoupe sur
$Y$ un ${\bf g}^1_2$ qui d\'etermine exactement le rev\^etement double
de $Y$ au dessus de $\pu$ d\'efini au paragraphe pr\'ec\'edent.

\textit{(\i\i)} Le quotient de $\p^3$ par une involution de la
famille $\mathfrak{I}$ est rationnel. Il est 
birationnel \`a $\p(H^0\I_X(4))$ 
: si $P$ est un point de $\p^3$ en dehors de $X$ alors on lui associe
l'ensemble des quartiques contenant $X$ et $P$ qui est un sous-espace de
codimension 1 de $H^0\I_X(4)$. R\'eciproquement, un sous-espace de
dimension 3 de $H^0\I_X(4)$ d\'etermine deux points points de $\p^3$
comme intersection r\'esiduelle de $X$.

\vs 0.4 cm

Nous allons maintenant d\'ecrire birationnellement cette famille
d'involutions modulo $PGL_4$ ou encore l'espace des modules des
courbes de degr\'e 9 et de genre 6 de $\p^3$.

Soit $X$ une courbe de degr\'e 9 et de genre 6 de $\p^3$. Nous lui
associons une cubique gauche $C$ de $\pd$ : l'ensemble des plans $H$
tels que $X\cap H$ est form\'e de 9 points situ\'es sur un pinceau de
cubiques est une cubique gauche $C$ de $\pd$. Ceci nous permet de
d\'efinir un morphisme $\P$ de $\H_{9,6}$ vers $\H_{3,0}$.

\vs 0.2 cm

\th{Proposition}{17} \textit{La fibre du morphisme $\P$ au dessus
d'une courbe $C\in\H_{3,0}$ est isomorphe \`a $\G(4,H^0\oo_{S^2C}(3))$
la vari\'et\'e des sous-espaces de dimension 4 de cubiques du plan
$S^2C$.}

\vs 0.2 cm

\dm :
\`A une cubique gauche $C$ nous pouvons associer le fibr\'e de
Schwarzenberger $F$ et r\'eciproquement (cf. proposition 5). Nous
avons vu qu'alors $\cv$ s'identifie \`a $S_3$ et $H^0F(1)$ \`a
$S^3S_2$ (c'est \`a dire les cubiques du plan $S^2C$).
Prenons un sous-espace vectoriel de dimension 4 de $H^0F(1)$. Nous
pouvons former la fl\`eche suivante : ${\check F}(-5)\fl \oo_{\p^3}(-4)^4$
dont le conoyau est l'id\'eal $\I_X$ d'une courbe de degr\'e 9 et
de genre 6.

R\'eciproquement, si $X$ est une courbe de degr\'e 9 et de genre 6 de
$\p^3$ nous avons la r\'esolution :
$$0\fl\op(-7)^2\fl\op(-6)^5\fl\op(-4)^4\fl\I_X$$
Nous pouvons alors consid\'erer le faisceau ${\check F}(-4)$ conoyau
de la premi\`ere fl\`eche. C'est un fibr\'e de Schwarzenberger (du
type de la proposition 5) associ\'e \`a la cubique gauche $C$ de
$\pd$ 
d\'efinie par 
l'ensemble des plans $H$ tels que $X\cap H$ est form\'e de 9 points
situ\'e sur un pinceau de cubiques.

\vs 0.4 cm

\th{Remarque}{10} Les transformations cubo-cubiques associ\'ees aux
instantons sont involutives. Elles d\'efinissent donc \'egalement des
involutions birationnelles de $\p^3$. Consid\'erons un r\'eseau $R$ de
quadriques. Soit $P$ un point de $\p^3$, son image par l'involution
est le point $P'$ orthogonal \`a $P$ pour toutes les formes
quadratiques du reseau $R$. Si $L$ est la droite qui joint $P$ \`a
$P'$, alors le reseau $R$ se restreint sur $L$ en une famille de
dimension 2 de formes quadratiques (sinon il n'existe pas de point
$P'$ sur $L$ qui est orthogonal \`a $P$ pour toutes les formes
quadratiques. Ceci signifie que la droite $L$ est contenue dans au
moins une quadrique du r\'eseau.

R\'eciproquement, si $L$ une droite contenue dans une quadrique du
r\'eseau, la restriction de $R$ \`a $L$ d\'efinit un pinceau de formes
quadratiques. L'othogonal de ce pinceau est une forme quadratique qui
a deux points isotropes. Ils sont reli\'es par l'involution.

Nous avons donc montr\'e que \textit{le quotient de $\p^3$ par l'involution
associ\'ee \`a $R$ est la vari\'et\'e des droites contenue dans une
des quadriques de $R$. C'est un complexe cubique de droites.}

On peut facilement v\'erifier que ce complexe cubique de droites est
un fibr\'e en coniques (toutes non singuli\`eres) au dessus d'une
surface de Del Pezzo de degr\'e 2 (rev\^etement double du plan
$\p(R)$ ramifi\'e au dessus de la quartique correspondant aux
c\^ones).  Ce fibr\'e en coniques ne semble pas avoir de section, je
ne sais pas si il est rationnel ou seulement unirationnel.


%



\begin{thebibliography}{99}




\bibitem[Ba]{Ba1} \textit{Wolf Barth} : Moduli of vector bundles on
the projective plane, Invent. Math. 42 (1977).


\noi
\bibitem[BB]{BB} \textit{Lionel Bayle et Arnaud Beauville} :
Birational involutions of $\p_2$. Kodaira's issue, Asian
J. Math. 4 (2000), no. 1.


\noi
\bibitem[Be]{Be} \textit{Eugenio Bertini} : Ricerche sulle
transformazioni univoche involutorie nel piano. Annali di Mat. 8,
224-286 (1877).


\noi
\bibitem[ES1]{ES} \textit{Geir Ellingsrud et Stein Arild Str\o mme} :
Stable rank-$2$ vector bundles on $\p^3$ with $c_1=0$ and $c_2=3$,
Math. Ann. 255 (1981). 


\noi
\bibitem[ES2]{ESGIT} \textit{Geir Ellingsrud et Stein Arild Str\o mme}
: On the Chow ring of a geometric quotient. Ann. of Math. (2) 130
(1989), no. 1.


\noi
\bibitem[FH]{FH} \textit{William Fulton et Joe Harris} :
Representation Theory, GTM Springer Verlag 129 (1991).


\noindent
\bibitem[GP]{GP} \textit{Laurent Gruson et Christian Peskine} :
Courbes de l'espace projectif : vari\'et\'es de s\'ecantes,
Enumerative Geometry and Classical Algebraic Geometry. Boston, Basel,
Stuttgart : Birkh\"auser (1982).
 

\noindent
\bibitem[GS]{GS} \textit{Laurent Gruson et Mohamed Skiti} :
3-instantons et r\'eseaux de quadriques, Math. Ann. 298 (1994).
 

\noindent
\bibitem[Ha1]{Ha1} \textit{Robin Hartshorne} : Algebraic geometry, GTM
Springer Verlag 52 (1977).


\noindent
\bibitem[Ha2]{Ha2} \textit{Robin Hartshorne} : Stable Reflexives
sheaves, Math. Ann. 254 (1980).
 
\noi
\bibitem[NT]{NT} \textit{Mudumbai S. Narasimhan et G\" unther
Trautmann} :  Compactification of ${\bf M}_{\p^3}(0,2)$ and Poncelet
pairs of conics. Pacific J. Math. 145 (1990), no. 2.


%


\noindent
\bibitem[OSS]{OSS} \textit{Christian Okonek, Michael Schneider et
Heinz Spindler } : Vector bundles on complex projective spaces, Basel,
Boston, Stuttgart : Birkh\"auser (1980).
 

\noindent
\bibitem[P1]{P1} \textit{Nicolas Perrin} : Une composante du bord des
instantons de degr\'e $3$, CRAS s\'erie I 330 (3) (2000).


\noi
\bibitem[P2]{P2} \textit{Nicolas Perrin} : Lieu singulier des surfaces
rationnelles. Pr\'ep. math.AG/0101083 (2001).


\noindent
\bibitem[P3]{P3} \textit{Nicolas Perrin} : Courbes rationnelles sur les
vari\'et\'es homog\`enes et une d\'esingularisation plus fine des
vari\'et\'es de Schubert. Pr\'ep. math.AG/0003199 (2000).




\noi
\bibitem[R]{R} \textit{Prabhakar A. Rao} : A family of vector bundles
on $\p^3$, LNM 1266 (1985).


\noi
\bibitem[Ro]{Ro} \textit{Maxwell Rosenlicht} : A remark on quotient
spaces, An. Acad. Brasil. Ci. 35 (1963).

\noindent 
\bibitem[RS]{RS} \textit{Leonard Roth et John G. Semple} :
Introduction to algebraic geometry, Oxford Science Publications
(1949).


\noi
\bibitem[S]{S} \textit{Rolf L. E. Schwarzenberger} : Vector bundles on
the projective plane, Proc. London Math. Soc. (3) 11 (1961). 


\noi
\bibitem[Sk]{Sk} \textit{Mohamed Skiti} : Espace de module des
fibr\'es vectoriels et groupe de Picard, en pr\'eparation.


\noindent
\bibitem[V]{V} \textit{Jean Valles} : Nombre maximal d'hyperplans
instables pour un fibr\'e de Steiner, Math. Zeit. 233 (2000).

\end{thebibliography}
\end{document}